\documentclass{article}
\usepackage[utf8]{inputenc}
\usepackage{amsmath,amsthm}
\usepackage{subcaption}
\usepackage[english]{babel}
\usepackage{graphicx,xcolor}
\usepackage[bottom,flushmargin,hang,multiple]{footmisc}
\usepackage{etoolbox}
\usepackage{url}
\usepackage{breakurl}
\usepackage[breaklinks]{hyperref}

\newtheorem{lemma}{Lemma}
\newtheorem{proposition}{Proposition}

\title{COVID-19 adaptive humoral immunity models: 
weakly neutralizing versus antibody-disease enhancement scenarios}
\author{
Antoine Danchin\footnote{The order of the authors is alphabetic}\footnote{Institut Cochin, Paris, France. \texttt{Antoine.Danchin@normalesup.org}}, 
\and 
Oriane Pagani-Azizi\footnote{CEREMADE - ESPCI Paris PSL, Paris, France. \texttt{Oriane.Pagani-Azizi@espci.fr}}, 
\and 
Gabriel Turinici\footnote{Université Paris Dauphine - PSL Research University, Paris, France. \texttt{Gabriel.Turinici@dauphine.fr}}, 
\and 
Ghozlane Yahiaoui\footnote{Mathematical Institute, University of Oxford, Oxford, United Kingdom. \texttt{Ghozlane.Yahiaoui@maths.ox.ac.uk}} 
}
\date{Jan 2022}
\begin{document}
\maketitle

\begin{abstract}
    The interplay between the virus, infected cells and the immune responses to SARS-CoV-2 is still under debate.
 Extending the basic model of viral dynamics we propose here a formal approach to describe the neutralizing versus weakly (or non-)neutralizing scenarios and compare with the possible effects of antibody-dependent enhancement (ADE). The theoretical model is consistent with data available from the literature;  
 we show that
  weakly neutralizing antibodies or ADE can both give rise to either final virus clearance or
   disease progression,  
 but the immuno-dynamic is different in each case.
 Given  that a significant part of the world population is already naturally immunized or vaccinated, we 
also discuss the implications on secondary infections 
 infections following vaccination or in presence of immune system dysfunctions.
\end{abstract}
%keywords: virus model
\section{Background}

SARS-CoV-2 is a new virus from the coronavirus family, responsible for the ongoing COVID-19 pandemic. To date, there are more than $300$ million cases  and  over five million deaths worldwide \cite{JHU}.  SARS-CoV-2 is the third betacoronavirus to severely infect humans appearing in the last 20 years, after SARS-CoV-1 and MERS-CoV. This motivates a growing need for efficient drugs and/or vaccines, not only for the time being but also in anticipation of a future coronavirus resurgence.

However, initial promising successes of antiviral treatments raised also the possibility of negative side-effects. On the vaccine front, an auto-immune disease (leading to temporary suspension of clinical trials) appeared during the AstraZeneca vaccine trial (Sept 9th 2020); altogether this context demonstrated the importance of understanding qualitatively and quantitatively the immune response to primary infection and also to challenges (vaccines belong to both categories). In particular, relevant mathematical models of the immune dynamics can be of interest to understand and predict the complicated behavior often observed.

We focus here on adaptive humoral immunity (the antibody-mediated immunity) and refer to future works  for an extension to the cellular and/or innate immune system. 

For clinical reasons and also for the understanding of those studying vaccines, antibody responses are of paramount importance. However, the neutralization capacities of antibodies is still under discussion, especially since weak or non-neutralizing antibodies can promote  infection through a process called antibody-dependent enhancement (hereafter abbreviated 'ADE') \cite{ADE1,ADE2,ADE3,ADE4} (see appendix \ref{sec:motivation_antibody} for details); thus the 
antibody neutralization capacity and ADE level are important ingredients of the model.
Furthermore, such behavior may  be accentuated by a challenge (secondary infection or infection following some immune system event or dysfunction), cf. also appendix \ref{sec:reinfection}.
Accordingly,
we investigated here both
primary and secondary COVID-19 infections. 

To summarize, 
we propose a mathematical model of the immune response and virus dynamics that includes the possibility of weakly neutralizing antibodies and / or ADE. {and analyze its implications.
	%Such situations can be accentuated by reinfection with a variant or for infections following vaccination. 
	At the time of writing the final draft (January 2022) a significant part of the world population is either vaccinated or naturally immunized and the consequences of  reinfection events are a major source of uncertainty concerning pandemic evolution.
This naturally calls for scientific investigation.}

\section{Methods}
\subsection{Mathematical model}

We present below the viral and immune response model. 
{It is a compartmental model similar to those used to describe the epidemic propagation, see \cite{kermack27,diekmann2000mathematical,MR1814049,ng_double_2003} for a general introduction and \cite{faraz_dynamic_2020,drozdzal_update_2021,liu_community_2020,danchin_immunity_2021,dolbeault_jean_heterogeneous_2020,biology10010010} for COVID-19 specific works.}

The viral-host interaction (excluding the immune response) is called the basic model of virus dynamics. It has been extensively validated both theoretically and experimentally, see \cite[eq (3.1) page 18]{nowak_virus_2000_book}, \cite[eqns. (2.3)-(2.4) page 26]{wodarz_killer_2007_book} and references therein. See also \cite{immo_revievs_math07,math_immunology_review16,eftimie_mathematical_immu_2016} for general overviews of mathematical immunology.

\begin{figure}[htbp!]
    \centering
    \includegraphics[width=0.90\textwidth]{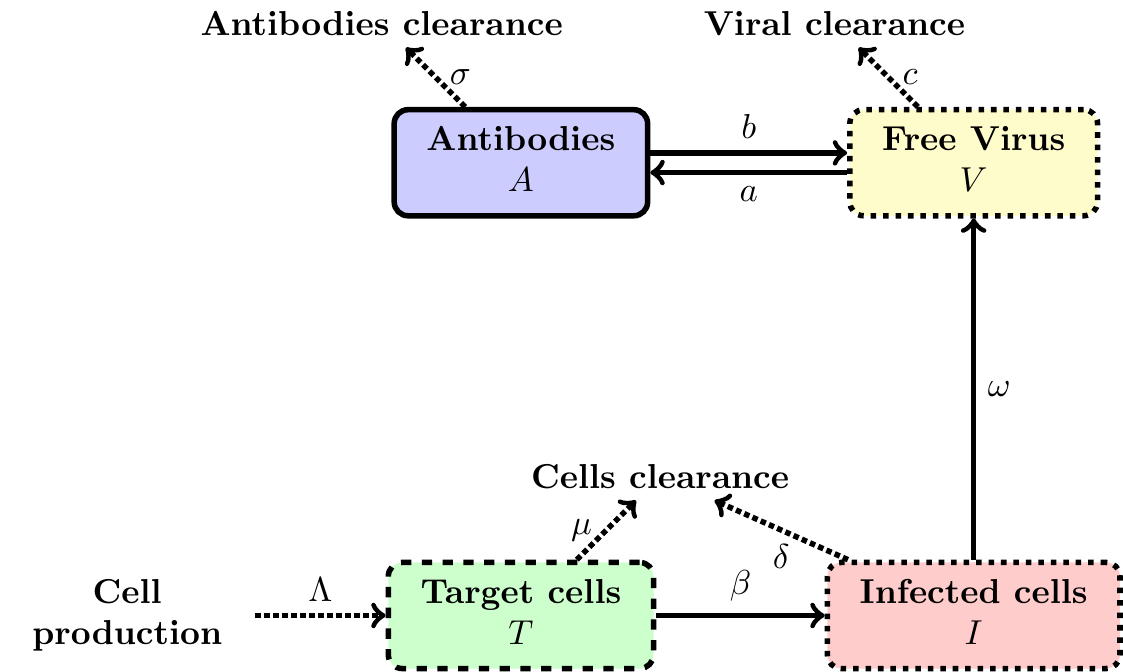}
    \caption{Graphical illustration of the flow in the model \eqref{eq:fullmodelT}-\eqref{eq:fullmodelbeta}.}
    \label{fig:tiva_model}
\end{figure}

The model involves several classes: that of the target cells, denoted $T$, the infected cells, denoted $I$, the free virus denoted $V$ and the antibodies denoted $A$. The model is illustrated in figure     \ref{fig:tiva_model}.

Target cells $T$, which in our case are the epithelial cells with ACE2 receptors located, for instance in the respiratory tracts including lungs, nasal and trachea/bronchial tissues, are produced at a rate $\Lambda$ and die at rate $\mu$. The parameters $\Lambda$ and $\mu$ define tissue dynamics in the absence of infection, see also section \ref{sec:noinfectionmodel}. When these susceptible cells meet free virus particles $V$, they become infected at a rate $\beta_0$. Furthermore, target cells can also become infected via ADE if virus entry is mediated by antibodies. The parameter $\beta_1$ represents the rate of ADE infection route which is the result of a three-species interaction: $T$, $A$ and $V$. 

Infected (initially target) cells, denoted $I$, die at a rate $\delta$. Note that this death rate will often be larger than the death rate of uninfected cells because viruses cause cell damage and cell death, \cite{wodarz_killer_2007_book,nowak_virus_2000_book}. Infected cells produce new virus particles at a rate $\omega$, and the free virus particles which have been released from infected cells decay at a rate $c$ called the clearance rate. 

Free virions are neutralized by antibodies A, which can block virus entry into cells but also facilitate phagocytosis, at a rate $b$. Finally, the antibodies can be stimulated by the free virus with a production rate $a$ while  declining at a rate of $\sigma$ (see for instance \cite[eq. (9.4) p.126]{wodarz_killer_2007_book}). {Note that alternative proposals for the antibody dynamics exist, see e.g. Andr\'e and Gandon \cite{andre_vaccination_2006} who 
	assume that
	immune response, once started, grows at a constant rate while 
	Pawlek et al. \cite{pawelek_within-host_2016} design a more complex model that takes into account the macrophage activation}.
The complete model reads (all constants are positive):
\begin{eqnarray}
& \ & 
dT/dt = \Lambda - \mu T -\beta(A) V T
\label{eq:fullmodelT}
\\ & \ & 
dI/dt =  \beta(A)V T - \delta I  \label{eq:fullmodelI}
\\ & \ & 
dV / dt = \omega I - c V - b A V  \label{eq:fullmodelV}
\\ & \ & 
dA/dt = a V A - \sigma A \label{eq:fullmodelA}
\\ & \ & 
\beta(A)=\beta_0+\beta_1 A. \label{eq:fullmodelbeta}
\end{eqnarray}

Several hypotheses in this model need to be further documented. The first one in that all infected cells including ADE infected cells support viral replication and can produce virus. However, to date, it is still unclear whether ADE infected cells can support viral replication in vivo, \cite{ADE3}, \cite{ADE4}. Here we choose not to distinguish between virus productive and non productive infected cells to keep the model simple
{(see however the comments in appendix \ref{sec:extended_latent_model}).} For the same reason, we do not discriminate between neutralizing, weakly neutralizing or non-neutralizing antibodies but consider all as members of the same class, the antibodies neutralizing capacity will therefore be the average of the neutralizing {power} and the average is described by the parameter $b$; on the other hand the ADE magnitude will be monitored by parameter $\beta_1$. These parameters are the most important part of the immune response and the object of our study.

{
\section{Stability of equilibria and further considerations} \label{sec:stability_main}

We operate under the assumptions that all parameters are positive and furthermore the following
two assumptions hold
(see Appendix \ref{sec:theory} for details):

\begin{eqnarray}
 & \ &   \textrm{Assumption 1:}	\ \ \ \ \ 
\delta > \mu.
	\label{eq:condition1} 
\\  & \ &
\textrm{Assumption 2:}	\ \ \ \ \ 
	\left( R_0 - 1 \right) \frac{\mu}{\beta_0}    > \frac{\sigma}{a}. 
	\label{eq:hypothesisV}
\end{eqnarray}
where
we define as in \cite[eq. (6.2) page 53]{nowak_virus_2000_book}):
\begin{equation}
	R_0 = \frac{\beta_0 \omega \Lambda}{c \delta \mu}.
	\label{eq:defR0}
\end{equation}
Note that \eqref{eq:hypothesisV} implies in particular $R_0 >1$ which is a standard condition for such models. We will further denote

\begin{equation}
V^{is}:= \left( R_0 - 1 \right) \frac{\mu}{\beta_0}, \ 
	V^{t} := \sigma/a.
	\label{eq:defVisVstar}
\end{equation}

\subsection{Stability of the equilibrium without ADE}
With these definitions we can give the main theoretical properties of the model depending on the
presence or not of the ADE term.

\begin{proposition}
The model
 \eqref{eq:fullmodelT}-\eqref{eq:fullmodelA} without ADE i.e.,
$\beta(A)=\beta_0$ (that is $\beta_1=0$) 
 has a single stable equilibrium given by:
\begin{equation}
T= T^{is}:=\frac{\Lambda}{\mu+\beta_0 V^t}, \ 
I= \frac{\beta_0 \Lambda V^t}{\delta(\mu+ \beta_0 V^t)}, \
V=V^{t}, \ 
A=\frac{c (V^{is}-V^t)}{\beta_0 b (\mu+ \beta_0 V^t)}.
\label{eq:equilibriumbeta1nul}
\end{equation}
\label{prop:stabilitynoade}
\end{proposition}
\begin{proof}
The proof of the stability of the equilibrium \eqref{eq:equilibriumbeta1nul}
is technical and is given in full detail in the appendix \ref{sec:modelvirusnoade}.
\end{proof}

\subsection{Stability of the equilibrium with ADE}

We investigate now the full model having a non-null ADE term $\beta_1 >0$.

\begin{proposition}
	The model
	\eqref{eq:fullmodelT}-\eqref{eq:fullmodelA} 
	has three equilibria:
\begin{enumerate}
	\item the trivial equilibrium $T = T^*=\Lambda/\mu$, $V=I=A=0$ which is unstable;
\item the immmunosuppression equilibrium, also unstable, given by~:
\begin{equation}
	T= T^{is}=\frac{\delta c}{\beta_0 \mu}, \ V=V^{is}:= \left( R_0 - 1 \right) \frac{\mu}{\beta_0}, \ I=I^{is}:=\left( R_0 - 1 \right) \frac{c \mu}{\omega \beta_0}.
	\label{eq:immunosuppression_eq}
\end{equation}
\item and a third equilibrium characterized as follows:
\begin{itemize}
	\item  the antibody level $A^f$ is the unique positive solution of the following 
	second order equation in the unknown $A$:
	\begin{equation}
		\omega \beta(A) \Lambda = \delta (c+bA) (\mu + \beta(A) V^t);
		\label{eq:equationAordre2}
	\end{equation}
	\item the other quantities are~:
	\begin{equation}
		T^f = \frac{\delta (c+bA^f)}{\omega \beta(A^f)}, \ 
		I^f = \frac{V (c+bA^f)}{\omega}, \ V=V^t=\frac{\sigma}{a}.
		\label{eq:equilibriumbeta1positive}
	\end{equation}

The following affirmations hold true concerning this third equilibrium
\begin{enumerate}
	\item when $\beta_1$ is small enough the equilibrium is stable;
	\label{item:prop_beta1_small}
	\item when $\beta_1$ is large enough the equilibrium is stable;
	\label{item:prop_beta1_large}
	\item however there exist choices of parameters (in particular values of $\beta_1$) for which this equilibrium is unstable. \label{item:prop_counter_example}
\end{enumerate}
		
\end{itemize}
\end{enumerate}

	\label{prop:stability_full}
\end{proposition}
\begin{proof}
	The proof is presented in appendix \ref{sec:modelvirus_full}.
\end{proof}
}

\subsection{Dynamical aspects} \label{sec:dynamicalaspects}

The equilibrium analysis in the previous sections does not yet tell the full story of the evolution of the system \eqref{eq:fullmodelT}-\eqref{eq:fullmodelA}.
Depending on the parameters, a common behavior is the following: 
initially $A$ will increase as response to $V$ being above threshold $V^t$; the increase of $A$ will drive both 
$I$ and $V$ to zero. Such a dynamics is stable over a long period and in practice $I$ and $V$ will keep small values for a time long enough to ensure virus clearance (when $V$ is small enough, due to the random nature of the events, $V$ will disappear).

Taking $I$ and $V$ to be constant equal to zero,  the new evolution is:
\begin{eqnarray}
	& \ & dT/dt = \Lambda - \mu T \\
	& \ & dA/dt =  - \sigma A.
\end{eqnarray}
Note that equations for $I$ and $V$ are missing because if the initial states are $V(0)=I(0)=0$ then $V(t)=I(t)$ for all $t\ge 0$. 
This evolution drives $T$ to $\Lambda/\mu$ and $A$ to zero.
If however during the slow decay of $A$ a challenge is presented in the form of a virus load $V>\sigma/a$ a new infection will start and $V$ and $I$ will rise again.

In conclusion, the stable equilibrium \eqref{eq:equationAordre2}-\eqref{eq:equilibriumbeta1positive}
is not necessarily reached in practice. The precise dynamics depends crucially on the parameters $b$ and $\beta_1$, see main text for details.

%Condition for continuous decrease of $I,V$ is 
%$\delta(c+bA) \ge \omega \beta_0 T$, realized in practice because $A$ is large.

%When $\beta_1>0$ can this prevent going to $I,V=0$?
%Condition: $\delta(c+bA) \ge \omega \beta(A) T$. When $A$ large this amounts to
%$\delta b \ge \omega \beta_1 T$

\section{Results}

\subsection{Theoretical results}

We refer the reader to the section \ref{sec:stability_main} for the rigorous statements concerning the theoretical properties of the model \eqref{eq:fullmodelT}-\eqref{eq:fullmodelbeta}. 
{Several situations may occur, but in summary the absence of ADE (i.e., $\beta_1=0$)
insures stable equilibrium while intermediate $\beta_1$ values (neither too small not too large) 
may provide examples of unstable equilibria; moreover,}
%We analyzed the equilibria of the model but the main conclusion is that 
stochastic events prevent the stable equilibrium state to be reached in practice, cf. section \ref{sec:dynamicalaspects}. The parameters $b$ and $\beta_1$ are shown to be the most important for the viral-host-antibody dynamics. 

{

\subsection{Empirical results: initial infection}
\label{sec:initial_infection}

Taking into account the available data from the literature and the methodology in appendix \ref{sec:choiceparams} 
%we used as baseline the parameters in table \ref{table:parameters}.
{
	we run a numerical procedure to fit the model parameters to reproduce at best the viral load data in figure \ref{fig:typ-variations} (left) and obtained the values in table \ref{table:parameters}. 
	The numerical simulation for a primary infection corresponding to 
	these parameters %in table \ref{table:parameters} 
	is shown in figure \ref{first_infection}.}

\begin{figure}[htbp!]
	\includegraphics[width=0.4\textwidth]{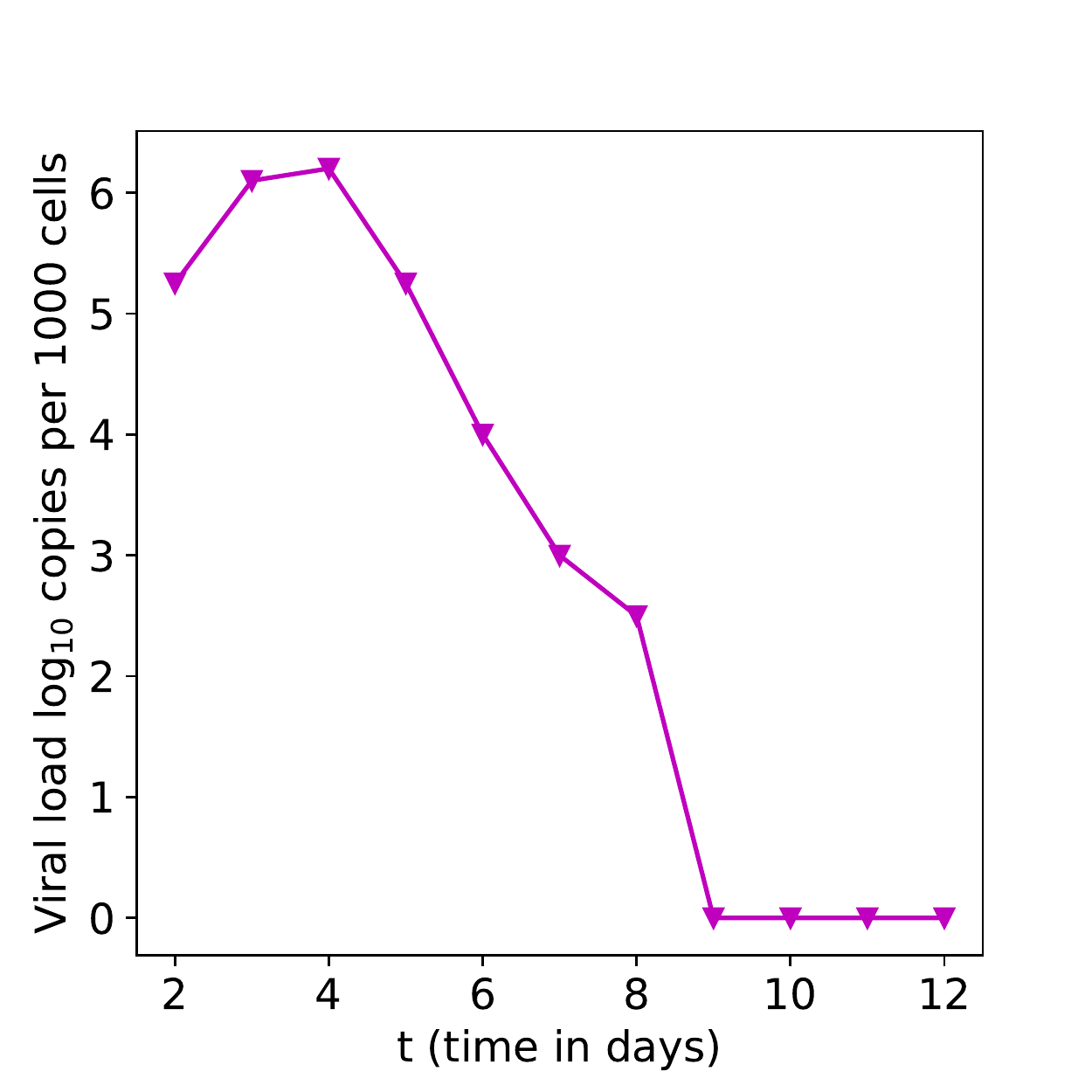}  
	\includegraphics[width=0.4\textwidth]{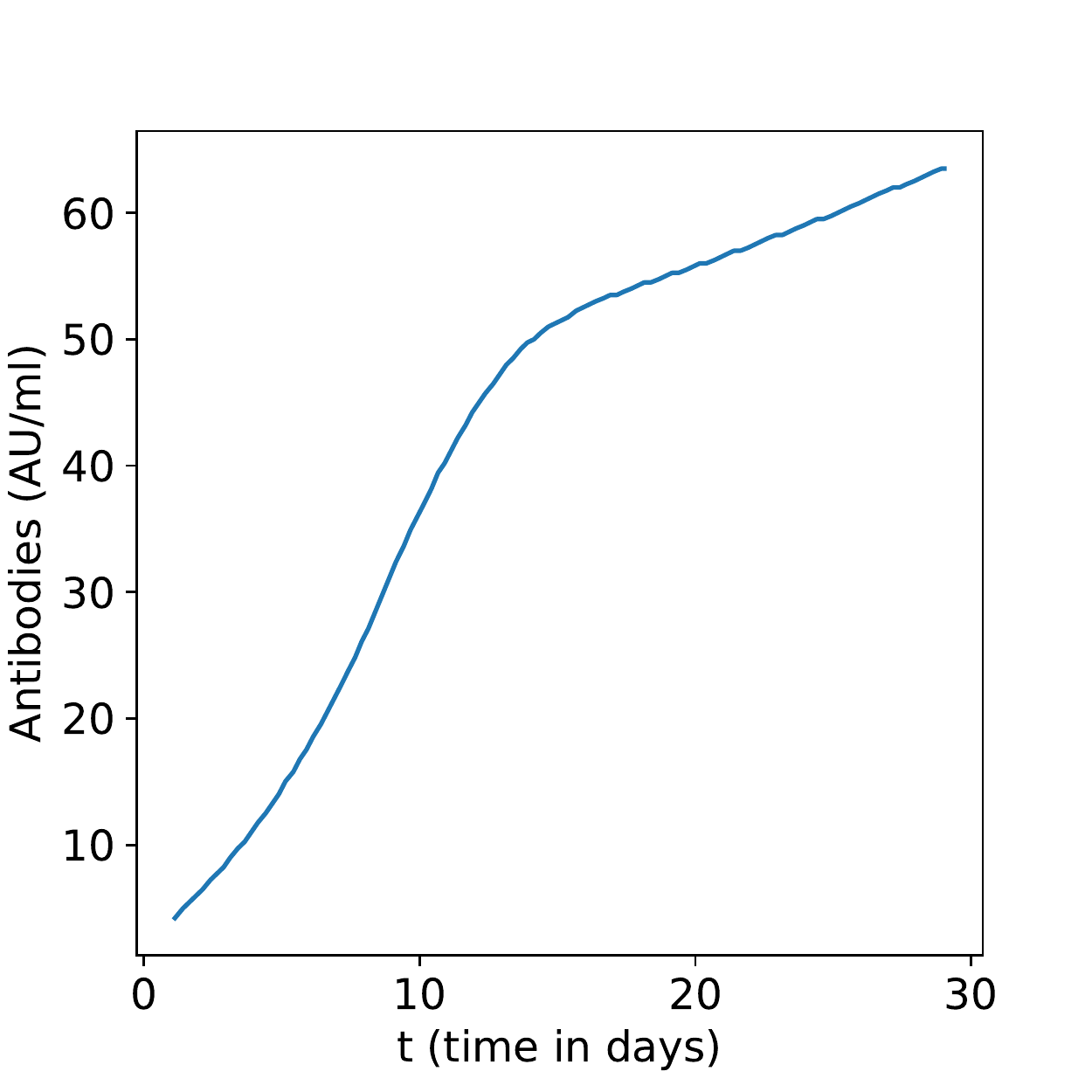}
	\caption{
		{\bf Left:}
		Clinically observed typical variation of SARS-CoV-2 viral load in nasopharyngeal swab normalised using cell quantification. Data taken from \cite[figure 3 page 703, patient 4]{Europe}. 	{\bf Right:}
		: Typical time variations for IgG. Data taken from \cite[figure 2 page 1085]{abs_figure}.
		{Note that the antibody data is a mean over several days and 
			corresponds to a different patient cohort.}
	}
	\label{fig:typ-variations}
\end{figure}

\begin{table}[hptb!]
	\begin{tabular}{|c|c|c|c|c|c|c|c|}
		\hline
		$\mu$ & $\Lambda$ & $\omega$ &  $\beta_0$& $\beta_1$ & $\delta$ & $c$  & $b$ \\
		
		\hline 
		$9.66$ & $9.66\times10^{6}$ & $59.74$ & $1.28\times10^{-6}$ & 0 &  $16.22$ & $1.45$ & $0.52$  \\
		\hline 
	\end{tabular}
	
	\noindent \begin{tabular}{|c|c|c|c|c|}
		\hline
		$a$&	$\sigma$ & $I(0)$ & $V(0)$ &  $A(0)$ \\
		\hline
		$9.15 \times 10^{-7}$ &$0.02$ & $372.11$ & $994.84$ &  $1.17$\\
		\hline
	\end{tabular}
	\caption{Baseline parameters use in numerical simulations of the model \eqref{eq:fullmodelT}-\eqref{eq:fullmodelA}.}
	\label{table:parameters}
\end{table}
}

\begin{figure}[htbp!]
    \centering
    \includegraphics[width=0.59\textwidth]{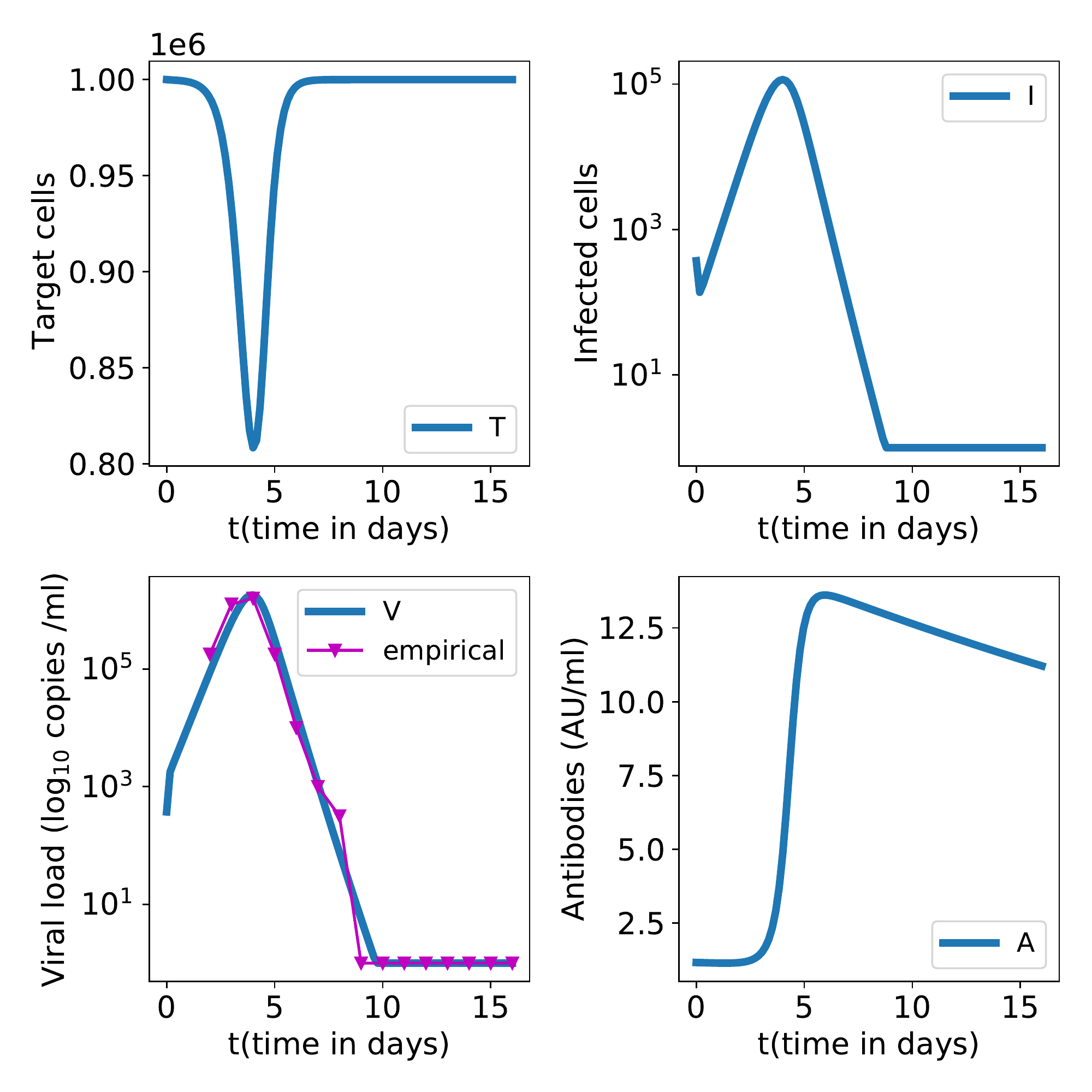}
    \caption{Numerical simulation of the first infection without ADE for model \eqref{eq:fullmodelT}-\eqref{eq:fullmodelA} and parameters in table \ref{table:parameters}.
{
A good fit for the viral load from figure \ref{fig:typ-variations} is obtained (data is truncated below the value $1$). 
On the contrary the fit is not as good for $A(t)$ because data 
does not correspond to the same patient (joint $V(t)/ A(t)$ data was not available).
}}
    \label{first_infection}
\end{figure}

There is a {$20\%$} fall of target cells which either become infected or naturally die. The viral load peaks around 4-5 days after symptoms onset at $1.8\times 10^{6}$ copies/ml. While SARS-CoV-1 viral load,  as MERS-CoV, peaked around 10 days after symptoms onset, most studies agree that SARS-CoV-2 viral load peaks sooner, around 5 days, \cite{Europe},\cite{upperresp}. Concerning antibodies, they increase sharply until week 2 then slower until a month after infection and start to decrease within 2-3 months \cite{sigmaest}, \cite{uk_Abs}. Qualitative agreement is observed with 
clinically observed variations variations of viral load and 
{partially with}
antibodies concentration depicted in figure \ref{fig:typ-variations} 
(see references in the figure).

{Note that, although we expect agreement 
	between $V(t)$ and 
 the viral load evolution in figure \ref{fig:typ-variations} (left) 
(which corresponds to a precise, real patient) the antibody
data from  \cite[figure 2 page 1085]{abs_figure} 
does not correspond to the same patient (data unavailable) but
is a mean value over several days and patients (not always the same). 
Each individual is likely to have his own
immuno-kinetic parameters: the parameters of the individual that may fit the $A(t)$ data from figure \ref{fig:typ-variations} (right) are not the same as the parameters that fit the data in left side of the same figure.}

The equilibrium state (\ref{eq:equilibriumbeta1nul}) when $\beta_1=0$ (no ADE present) is reached after $2$ years for all variables in figure \ref{first_infection}. However, viral load and infected cells reach a minimum {within several weeks} post-infection before increasing and oscillating toward equilibrium state (\ref{eq:equilibriumbeta1nul}) 
(simulations not shown here).
Therefore, if the virus load is very small close to the minimum, all other variables decrease towards 0 and the infection has vanished. The equilibrium state (\ref{eq:equilibriumbeta1nul}) is stable but not reached in practice as the patient is cured. 

\subsection{Empirical results: secondary infection, variants, vaccination}
\label{sec:secondary_infection}

We focus {on a scenario where the immuno-kinetic parameters such as the neutralizing efficacy ($b$) or the ADE parameter ($\beta_1$) change; the causes can be multiple: a primary infection with a different variant, vaccination, or some immune evolution (aging  being an example). In all cases we investigate the infection, called challenge, that takes place with a different set of $b$ of $\beta_1$ parameters than in table \ref{table:parameters}.}

\subsubsection{Variation of the neutralizing capacity $b$}
\label{sec:secondary_infection_b}

\begin{figure}[htbp!]
	\centering
	\includegraphics[width=0.59\textwidth]{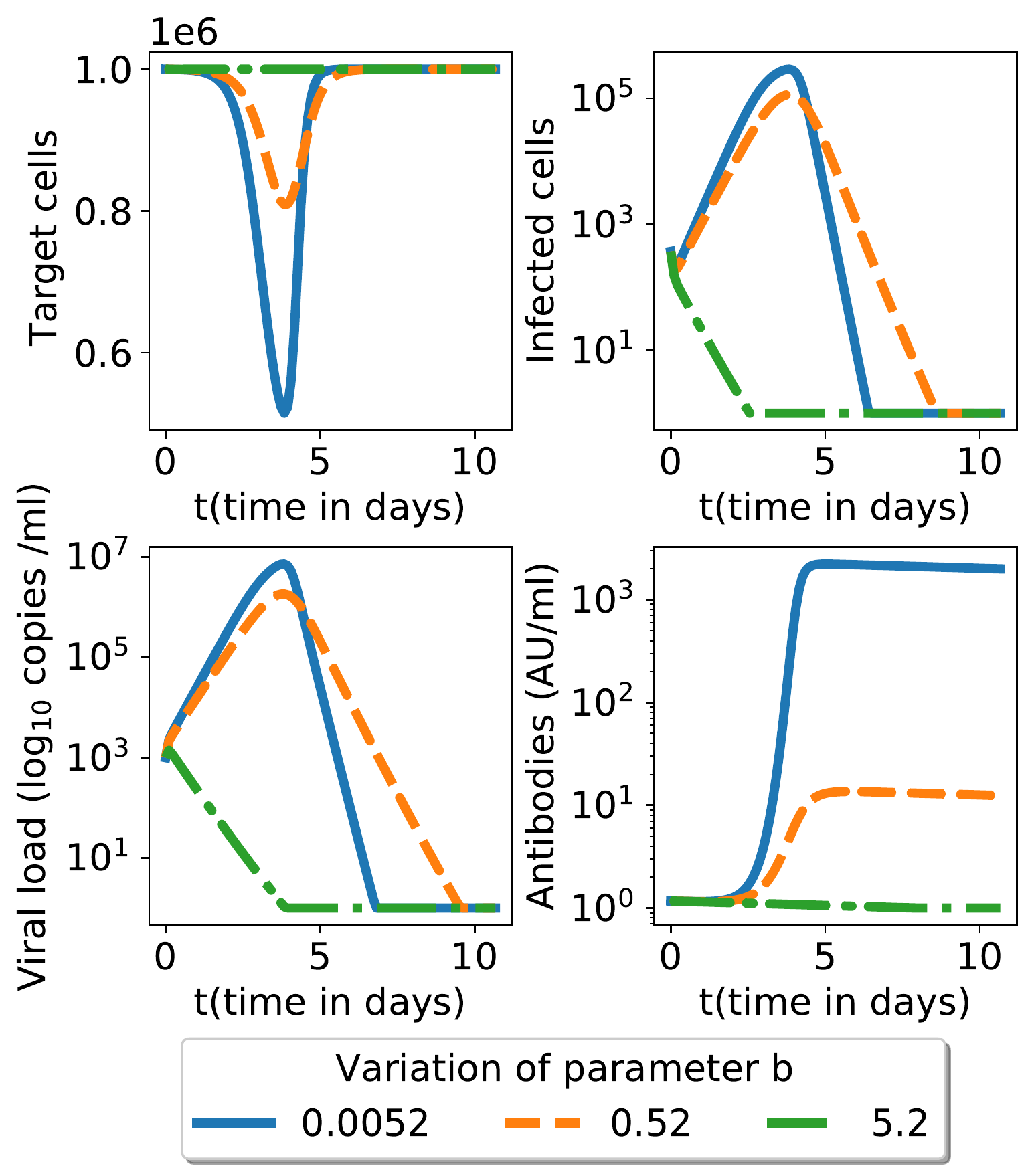}
	\caption{Numerical simulation of the model \eqref{eq:fullmodelT}-\eqref{eq:fullmodelA} and parameters in table \ref{table:parameters}.
		{Only the neutralizing capacity $b$ is changed around the nominal value $b=0.52$. Lower value of the neutralizing capacity $b$ (solid blue line) leads to lower target cell count and higher antibody levels.
See figure \ref{fig:variable_blong} for simulation over a longer time span.}}
	\label{fig:variable_b}
\end{figure}

When there is no ADE, decreasing the neutralizing capacity of antibodies 
{(parameter $b$)} leads on the one hand to a higher viral load peak but on the other hand to higher antibodies concentrations. The less neutralizing the more abundant antibodies are to compensate so that the infection is always solved.
{The simulations results are presented in figure \ref{fig:variable_b}.
Infection resolution is obtained with little target cell destruction for larger values of $b$. On the contrary, low values of 
$b$ will lead to significant increase of the antibody number and  simultaneous decay of target cells, both largely pejorative for the patient.

In the cases where the viral load reaches low values the infection stops before converging to the theoretical equilibrium.
}

\subsubsection{Presence of ADE ($\beta_1>0$)}
\label{sec:secondary_infection_ade}

\begin{figure}[htbp!]
    \centering
    \includegraphics[width=0.59\textwidth]{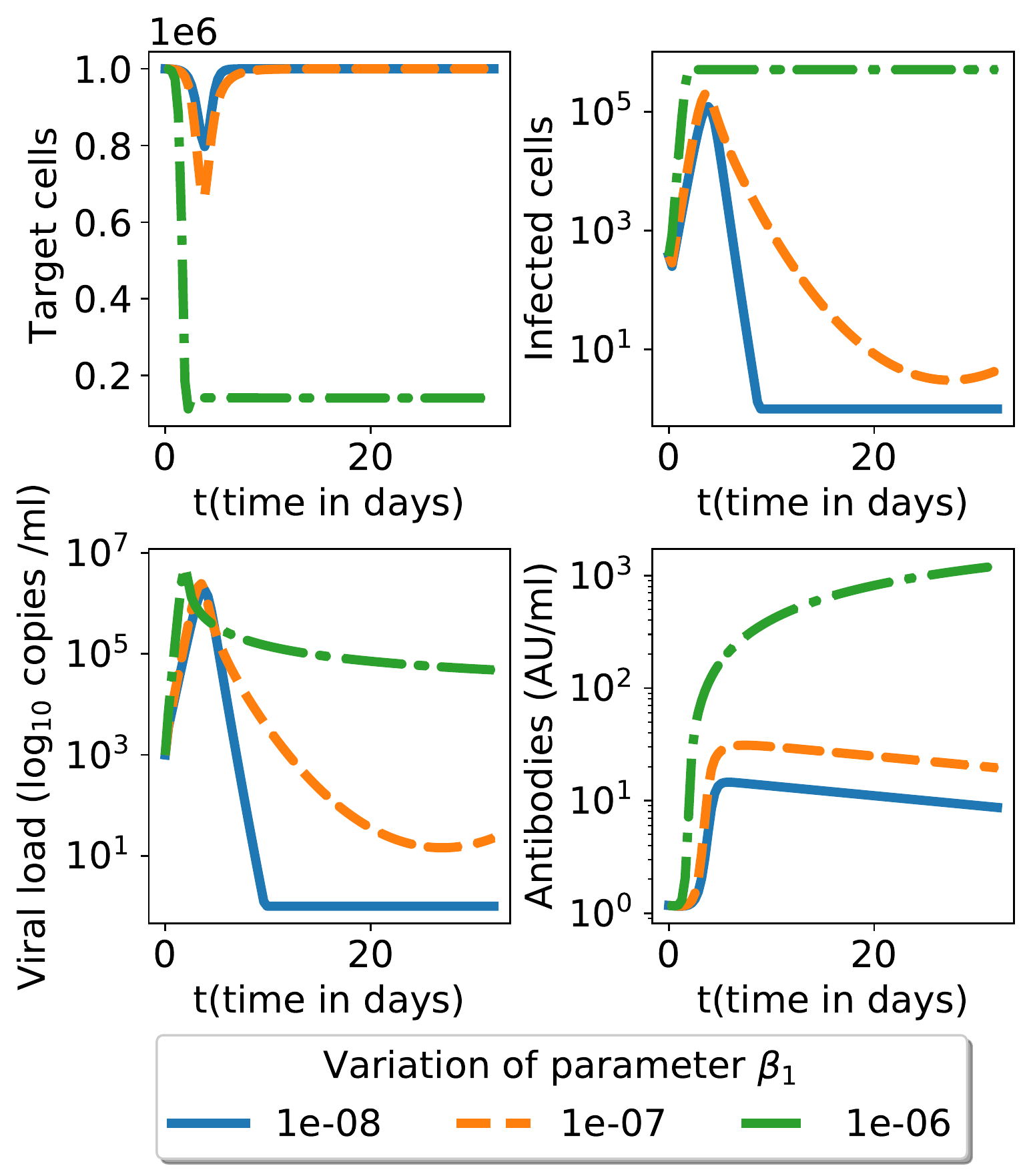}
    \caption{The effect of the ADE parameter $\beta_1$ (the model \eqref{eq:fullmodelT}-\eqref{eq:fullmodelA}). The secondary (challenge) infection has fixed neutralizing antibodies capacity $b$ but several
    	{possible} 
    	 $\beta_1$ (ADE) parameter values; for all other parameters  we use the nominal values given in table \ref{table:parameters}.
{See figure \ref{fig:VariationADElong} for simulation over a longer time span.     
}
}
    \label{fig:VariationADE}
\end{figure}

We investigated in
figure \ref{fig:VariationADE} 
the possibility of the ADE mechanism present ($\beta_1>0$), for a range of possible parameter 
$\beta_1$ values. 
We plot all variables upon challenge with the same neutralizing capacity for antibodies. 
%ADE can be triggered by several mechanisms, such as weakly neutralizing antibodies or sub-optimal concentration of neutralizing antibodies, which, for simplicity, are not distinguished in this model.
 A higher ADE parameter leads to more destroyed target cells, more infected cells, more viral load and more antibodies. However the antibodies concentration is restricted by an upper limit. There is a threshold effect : increasing $\beta_1$ does not increase significantly the antibody population (see figure \ref{fig:VariationADE} and compare with theoretical insights in the proof of point \ref{item:prop_beta1_large} of proposition \ref{prop:stability_full} in Appendix \ref{sec:theory}). Therefore a higher $\beta_1$ ADE parameter cannot be compensated by more antibodies. 
 %as a lower neutralizing capacity was already present for the first infection. 
  For example, unlike $\beta_1=10^{-8}$, if $\beta_1=10^{-6}$ the viral load directly stabilizes to its equilibrium state (\ref{eq:equilibriumbeta1positive}), without reaching a minimum close to 0 while oscillating (simulation not shown here). In this case, the infection wins (leading to respiratory function disruption and possibly patient death).
{Large values of $\beta_1$ lead to significant (possibly total) destruction of target cells.}

\section{Discussion}

We investigated the immune response
to SARS-CoV-2 infection and re-infection through a numerical model; the model 
 can also take into account the possible presence of ADE, either on first infection or to a challenge (secondary or re-infection with a different phenotype, after vaccination, etc.). 

As to date there is no clear evidence that ADE occurs in COVID-19 severe patients, we assume that ADE only happens upon challenge. %However, if solid evidence were given, our results for challenge still hold for the first infection.  

We started from a classic viral-host dynamics (\cite{wodarz_killer_2007_book},\cite{DENV/ZIKV}) that we modified by adding parameter $\beta_1$ to account for a possible ADE mechanism. In order to keep the model at its lowest complexity, we do not distinguish between ADE triggering antibodies and neutralizing antibodies.

We conducted a theoretical study of our system by computing equilibrium states and stability with and without ADE. We showed that stochastic events may {also play a role and} prevent the stable equilibrium state to be reached in practice; we identified the parameters $b$ (neutralization capacity) and $\beta_1$ (ADE presence) to be crucial for the dynamics of our system. 

Then, we calibrated our parameters values to match reference viral load from the literature
\cite{Europe,uk_Abs} {and obtained good results.}
% and antibodies concentration. 

{We then investigated a secondary infection (or infection following vaccination or other immune event) that can posess different immuno-dynamic parameters.} We saw that without ADE, the possible weak antibody neutralizing capacity was systematically compensated with higher concentrations of antibody leading to viral clearance. On the other side, adding ADE was not always associated with viral clearance {but possibly high target cell destruction.} Simulations and equilibrium analysis showed that antibody concentration had an upper limit which prevented higher ADE to be compensated by an unlimited antibody quantity. Therefore, ADE should be taken in consideration as a serious risk in 
disease understanding, treatment {and vaccine development and scheduling.}

On the other hand, we showed that the results are sensitive to the neutralizing antibody capacity (the $b$ parameter); note that a decrease of this parameter can occur in several situations, for instance due to immune function decay, due to the malfunctioning of the antibody immunodominance mechanism that ends up selecting too many weakly neutralizing antibodies or due to miscalibrated therapeutic interventions. Independent of the cause, such a decrease of the neutralizing capacity is susceptible to imply a substantial deterioration of the outcome.

{In summary, our results seem to support a picture where ADE presence is 
primary correlated with important target cell destruction 
while loss in neutralization capacity is correlated with both higher antibody count (leading to inflammation) and larger target cell destruction.}

\subsection{Limitations and future work}

As any other, our model contains of course several limitations. First, we considered all infected cells to support viral replication, including ADE-infected cells. Concerning SARS-CoV-2, the questions of ADE is still under debate, but for SARS-CoV-1 in vitro ADE evidence suggested abortive viral replication in ADE infected cells. Therefore, if we changed the model \eqref{eq:fullmodelT}-\eqref{eq:fullmodelA} to include this distinction, equilibrium state would change and ADE may be compensated. Similarly, we did not distinguish between former antibodies and novel antibodies secreted upon challenge. This would imply more parameters and change equilibrium levels but without inherently changing variables behavior. 
{The antibody dynamic model can also be changed to include e.g., constant antibody
	production following a threshold or more specific effects \cite{andre_vaccination_2006,pawelek_within-host_2016}.}
Regarding parameters, we did not have enough exploitable available data to train our model and fit parameters better. 
Finally, an unique model can hardly capture the extreme variability of COVID-19 clinical outcomes, see \cite{questions_nature}; some studies proposed that some of the variability come from genetics, see e.g., \cite{genetics} where genetic information from roughly 4,000 people from Italy and Spain  was  correlated to severity of COVID-19. This may lead to a variability of our model parameters in the form of random variables.
%Respiratory failure was associated with one of two particular gene variants : one lies in the region of the genome that determines ABO blood type, the other one is near several genes, including one that encodes a protein that interacts with the receptor the virus use
%s to enter human cells, and two others that encode molecules linked to immune response against pathogens. 

{
By now (January 2022), billions of people have been vaccinated, at
least by one injection, and more than half a billion have been infected. 
Fortunately, despite the
spread of the highly contagious omicron variant, it appears that
morbidity and mortality are declining. This implies that, for the time
being, the most disastrous consequences of the phenomena 
included in this model are not being observed. 
However, it must be emphasized that human
polymorphism, viral polymorphism and highly variable environmental
conditions, as well as the considerable variety of vaccination
protocols, mean that there may be isolates where ADE or the other immune
responses we have explored could be significant. It is therefore
particularly important to monitor variations in morbidity and mortality
around the world so that a rapid response can be implemented if there is
any local increase. Finally, the types of vaccines used are very
different. For those based on well-established technologies, we do not
foresee any consequences other than those discussed in this work, except
perhaps as a function of the vaccination protocols (time lag between
primary and booster injections). In contrast, the use of vaccines based
on indirect antigen production (adenovirus- or synthetic RNA-based
vaccines) requires specific encapsulation of the active ingredient in a
variety of capsules or cassettes. These containers can, by themselves,
be immunogenic. The consequence would be that after several
immunisations, patients would develop a response against the vaccine,
rendering it ineffective against the disease. We did not consider this
outcome in our work.
}

The more science will shed light on the full picture of SARS-CoV-2, the more our model can input complex and precise details. In the meantime, the main take-home message is that, with parameters consistent with the available clinical data,
the neutralizing capacity and ADE mechanisms may play an important immunological role into the primary and secondary infection outcomes.

\section*{Acknowledgements}

Ghozlane Yahiaoui is supported by the Engineering and Physical Sciences Research Council (EPSRC): CDT Grant Ref. EP/L015811/1.

%%\bibliographystyle{plain}
%\bibliographystyle{unsrt}
%\bibliography{references}
%biblatex style
%\printbibliography

\newpage
\appendix

\section{Additional motivations}

\subsection{Available evidence detailing the antibody response}
\label{sec:motivation_antibody}

For clinical reasons and also for the understanding of those studying vaccines, antibody responses are of paramount importance; SARS-CoV-2 specific antibodies are usually detected during the second week after illness onset (see \cite{AbResp7,AbResp5}) and remain active thereafter for an unknown time span (see 
however \cite{seow_longitudinal_2020} for 
recent information). Antibody responses are mainly directed against the RBD-spike and nucleocapsid proteins. However, the neutralization capacities of these specific antibodies is still under discussion, especially since weak or non-neutralizing antibodies can enhance infection through a process called antibody-dependent enhancement (hereafter abbreviated 'ADE') \cite{ADE1,ADE2,ADE3,ADE4}. This has been recently emphasized in the set up of clinical trials (see for example 
\cite{moderna}), in a general discussion of the prospects of vaccination \cite{defrancesco_whither_2020} and in a perspective accounting for the present situation in terms of SARS-CoV-2 vaccines, therapies and immunity \cite{lee_antibody-dependent_2020,poland_sars-cov-2_2020}

The present academic interpretation of the ADE is that it occurs through virus-antibody immunocomplexes that facilitate virus internalization in host cells that do not express virus receptor but Fc receptors. ADE is induced when the antibody-virus stoichiometry is below the threshold for neutralization, \cite{ADE1,ADE2}. As a consequence, neutralizing antibodies may enhance infection when their concentrations fall below a key occupancy threshold, and some poorly neutralizing antibodies may strongly increase infection over a wide dose–response range. ADE has been demonstrated in vitro for many viral infections, including that triggered by SARS-CoV-1 which was reported to infect in vitro human macrophages (see \cite{ADE3}) and human B cell lines via an ADE pathway, (see \cite{ADE4}). Moreover, Qidi Wang et al. reported that a specific spike protein epitope elicited antibodies which could enhance infection via ADE, while other epitopes induced neutralizing antibodies in non-human primates. Furthermore, the authors showed that a SARS-CoV-1 inactivated vaccine could induce ADE and lung pathology in experimental rhesus monkeys \cite{ADE5}. In contrast, Martial Jaume and co-authors showed that vaccine candidates which mediated in vitro ADE infection could still be neutralizing and protective in vivo on rhesus macaques, \cite{ADE4}. Moreover, in most cases, ADE infected cells do not support viral replication, \cite{ADE3,ADE4}. Instead, ADE may trigger cell apoptosis and promote tissue inflammation and injury with the release of pro-inflammatory cytokines from infected cells, \cite{ADE2,ADE6}. As a result, whether ADE actually happens in SARS-CoV-1 infected humans and is a factor of disease severity is still a debated research subject since no in vivo human evidence has been demonstrated yet (but this statement is very time-dependent given the present intense research on SARS-CoV-2). Note however, that SARS-CoV-1 infected patients who developed a higher and earlier antibody response were associated with worse clinical outcome. An early antibody response may be weakly neutralizing compared to a later one. As a consequence, a high concentration of those antibodies could lead to ADE and enhancement of infection.  

The question of ADE and the link between antibody dynamics and disease evolution is still unclear for COVID-19. J. Zhao et al. reported a strong positive correlation between disease severity and high antibody titers two weeks after illness onset. The antibody level is considered as a risk factor for severe evolution, independently from age, gender and comorbidities \cite{AbResp1}. In another study, Wenting Tan et al also came to the same conclusion: higher titers of anti-N IgM and anti-N IgG are observed for severe patients \cite{AbResp2}. Finally, Baoqing Sun et al observed that severe patients had higher levels N-IgG than S-IgG after the symptoms onset. However, according to the authors, whether N-specific antibodies can block virus infection is still open to question \cite{AbResp3}. The secretion of a high level of weakly neutralizing antibody supports the hypothesis of ADE for COVID-19 which can partially explain some clinical complications. In contrast, Mehul S. Suthar et al concluded that the appearance of high titer neutralizing antibody responses early after the infection was promising and may offer some degree of protection against re-infection \cite{AbResp4}. This result seems to be confirmed in a recent study in which SARS-CoV-2 infection induced protective immunity against re-exposure in nonhuman primates. However, rhesus macaques do not develop severe clinical complications as reported in human patients, suggesting that if rhesus macaques produce neutralizing antibodies, transposition of this observation to humans is still to be investigated \cite{rhesussars2}. Finally, a recent study on a recovered cohort of COVID-19 patients showed that elderly patients had significantly higher levels of antibodies than younger patients. However, severe and critical patients were excluded from the study because they received passive antibody treatment before sample collection. As a result, the authors could not directly evaluate the effect of antibodies on virus clearance or disease progression in COVID-19 patients \cite{AbResp6}. This suggests that if elderly patients tend to develop higher titers of antibodies, those may not be systematically associated with worse clinical evolution. What should rather be answered is whether disease severity is systematically associated with high antibody levels. 

On the other hand, the vaccine community is increasingly aware of this need (see discussion on ADE in \cite{lee_antibody-dependent_2020,defrancesco_whither_2020,poland_sars-cov-2_2020}) and studies along these lines are required.

Another motivation comes from the fact that 
the adaptive immune system response starts in about a week; on the other hand in many mild forms infection is resolved in around a week while on the contrary severe forms may at first start as mild and only then become severe; a simplistic view may indicate that the innate immune response is very efficient while the adaptive immune system response may be detrimental. In this case, everyone with a mild first infection (i.e., mostly dealt with by the innate immune system) will, upon re-infection, see an adaptive immune response rising faster (once the memory is in place, its response is faster than the innate immune response) and thus the detrimental effects could be visible for people previously having experienced mild forms, e.g., low age class individuals.

\subsection{Evidence on re-infection} \label{sec:reinfection}

The possible unfavorable outcomes of a secondary infection (challenge) following a primary SARS-CoV-2 infection were described in various situations (see for example
\cite{larson_case_2021}), but an increasing body of evidence highlights the Kawasaki-like syndrome as a possible negative outcome, see \cite{kawasaki20,courage_what_2020,Bergame,riphagen_hyperinflammatory_2020,viner_kawasaki-like_2020}.

An italian study \cite{Bergame} indicates that the immune response to SARS-Cov-2 is responsible for the appearance of a pediatric Kawasaki-like syndrome (Kawasaki-like disease or Multisystem Inflammatory Syndrome in Children MIS-C in the US). In this study, 8 to 10 children have been tested positive to IgG, IgM or both (the infection to SARS-CoV-2 preceded the development of the syndrome) and 2 only in PCR (the infection was simultaneous). SARS-CoV-2 infected children who developed the Kawasaki-like syndrome (KLS) were on average older and more severely hit than other children victims of the classical Kawasaki syndrome.

The same phenomena has been observed in the US and UK \cite{kawasaki20,riphagen_hyperinflammatory_2020,viner_kawasaki-like_2020}. Academic studies begin to investivage the interplay between COVID-19 and the MIS-C \cite{consiglio_immunology_2020}.

The causes of the development of the Kawasaki disease are still unknown. The best accepted hypothesis is that of an abnormal immune response that occurs as a result of the infection provoked by one of several pathogenic agents
%(facteurs environnementaux qui déclenchent l'apparition de maladies autoimmunes cf le powerpoint AID-ADE) 
for the genetically susceptible patients. The triggering pathogens have not yet been identified. 

To account for the peaks of the KLS cases following an infection with SARS-CoV-2, two hypothesis may be formulated: the antibodies produced by the children can induce the initiation of an autoimmune disease and  syndromes similar to the Kawasaki syndrome. The second hypothesis is an ADE-type mechanism.

\section{Choice of simulation parameters}
\label{sec:choiceparams}

{
Parameters' order of magnitude were derived from literature, see \cite{paramrefs} for $\mu$, \cite{Europe} for $\omega$,
clearance data from \cite{DENV/ZIKV,parametersSARS2,sigmaest,uk_Abs}. To obtain the precise values, 
we then fitted the model to the SARS-CoV-2 clinical data available in figure \ref{fig:typ-variations} and obtained the values in table \ref{table:parameters} (simulation results are shown in figure \ref{first_infection}).}

\section{Sensitivity with respect to parameters}
\label{sec:supplementary}

{We plot here a longer time evolution corresponding to figures \ref{fig:variable_b} and \ref{fig:VariationADE}. This allows to see the difference between initial dynamics and the long time equilibrium, cf. considerations in section \ref{sec:dynamicalaspects}.}

\begin{figure}[htbp!]
    \centering
    \includegraphics[width=0.59\textwidth]{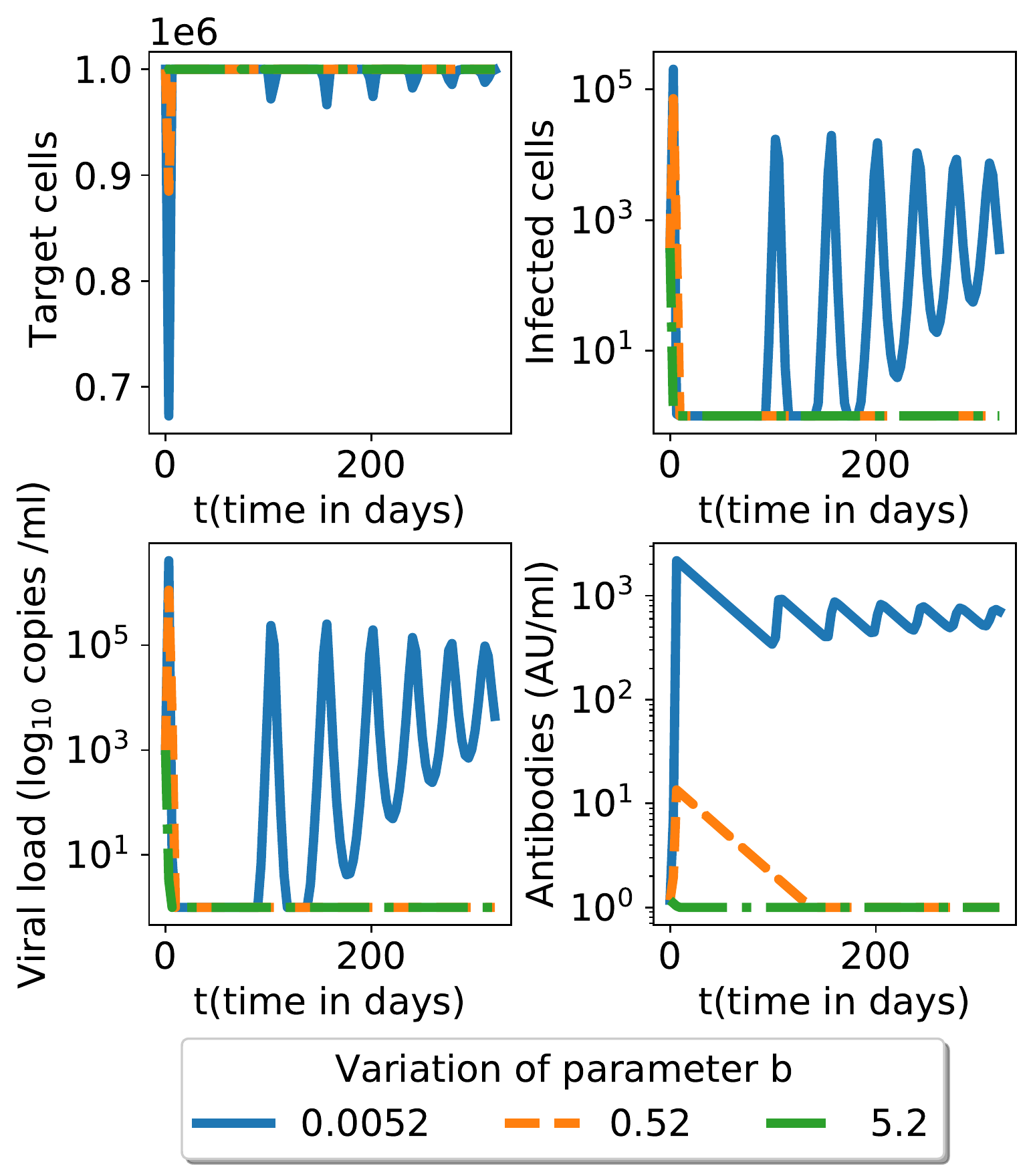}
    \caption{Simulation in figure \ref{fig:variable_b} for a longer time span.}
    \label{fig:variable_blong}
\end{figure}
\begin{figure}[htbp!]
	\centering
	\includegraphics[width=0.59\textwidth]{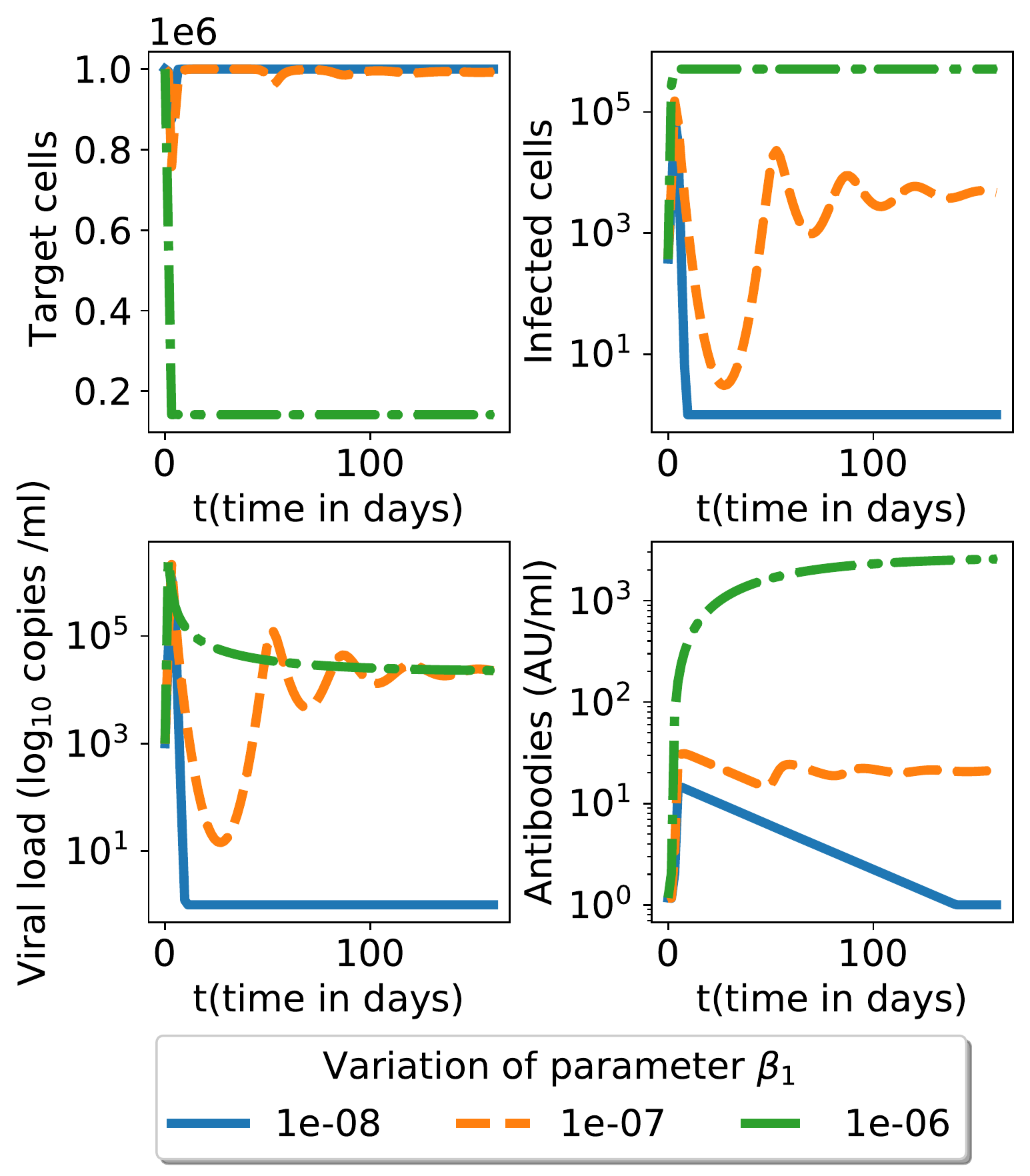}
	\caption{Simulation in \ref{fig:VariationADE} for a longer time span.}
	\label{fig:VariationADElong}
\end{figure}

\section{Mathematical properties of the model} \label{sec:theory}

We describe in an incremental way the mathematical properties of the main model \eqref{eq:fullmodelT}-\eqref{eq:fullmodelbeta}. {We take advantage of this description
	to illustrate the hypotheses \eqref{eq:condition1} and \eqref{eq:hypothesisV}. 
	The results in sections \ref{sec:noinfectionmodel} and \ref{sec:modelvirusnoade} are known, 
	see e.g.,
	\cite{perelson_dynamics_1993,nowak_virus_2000_book,smith_virus_2003,wodarz_killer_2007_book} while those in the main text (propositions \ref{prop:stabilitynoade}, \ref{prop:stability_full} and their proofs in this appendix) are, to the best of our knowledge, original.
}

\subsection{Model without a virus, nor immune response}
\label{sec:noinfectionmodel}

In absence of any infection the equations for the target cells are (see \cite{wodarz_killer_2007_book,nowak_virus_2000_book}):
\begin{eqnarray}
	\label{eq:Model_ADE_virusfree}
	& \ & 
	dT/dt = \Lambda - \mu T.
\end{eqnarray}

Since the Jacobian matrix at equilibrium (a $1\times 1$ matrix) is the constant $-\mu$ therefore the equilibrium is stable, in fact 
any initial data $T(0)$ will converge to the equilibrium  
\begin{equation}
	T^* = \Lambda/\mu.    \label{eq:equilibrium1d}
\end{equation} 

\subsection{Model with virus but no immune response}
We employ the basic model of virus dynamics, see \cite[eq (3.1) page 18]{nowak_virus_2000_book} and also \cite[eqns. (2.3)-(2.4) page 26]{wodarz_killer_2007_book} described by the equations:
\begin{eqnarray}
	& \ & 
	dT/dt = \Lambda - \mu T -\beta_0 V T, 
	\label{eq:Model_ADE_virus_no_aT}
	\\ & \ & 
	dI/dt =  \beta_0V T - \delta I, 
	\label{eq:Model_ADE_virus_no_aI}
	\\ & \ & 
	dV / dt = \omega I - c V. 
	\label{eq:Model_ADE_virus_no_aV}
\end{eqnarray}
The initial conditions are:
\begin{equation}
	T(0) = T^*=\Lambda/\mu, \ I(0) = 0, \ V(0) > 0,
\end{equation}
which express the fact that the initial state for $T$ is the stable equilibrium seen in section \ref{sec:noinfectionmodel}, there are initially no infected cells and the initial viral load is strictly positive.

It is natural to assume that the decay rate of infected cells is at least as large as the decay rate of healthy cells, i.e., assumption \eqref{eq:condition1}.

In this model, an infection is only possible if the basic reproduction ratio of the virus in the absence of immune response, defined in \eqref{eq:defR0}
is strictly super-unitary, that is 
\begin{equation}
	R_0 > 1. \label{eq:condition}
\end{equation}
Otherwise, that is if $R_0\le 1$, the initial viral load can only decrease. 
The model has two equilibria: 

- trivial equilibrium: $T = T^*=\Lambda/\mu$, $V=I=0$. The 
Jacobian matrix at equilibrium is $\begin{pmatrix} -\mu & 0 & -\beta_0 T^* \\ 0 & -\delta & \beta_0 T^* \\ 0 & \omega & -c \end{pmatrix}$. The eigenvalues of this matrix, under condition \eqref{eq:condition}, are all real but not all negative: one of them is $\lambda_1=- \mu$ but the product of the other two is $\delta c - \omega \beta_0 T^* \le 0$ thus at least one is positive. Therefore, {under assumption \eqref{eq:condition},} this critical point is not {a stable} equilibrium.

- the "immunosuppression" equilibrium 	\eqref{eq:immunosuppression_eq}

%For values in table \ref{table:parameters}, we obtain : $T^{is}=4.55\times10^{8}$, $V^{is}=2.18\times10^{8}$ and  $I^{is}=3.28\times10^{5}$.

The Jacobian matrix is  
$\begin{pmatrix} -\mu-\beta_0 V^{is} & 0 & -\beta_0 T^{is} \\ \beta_0 V^{is} & -\delta & \beta_0 T^{is} \\ 0 & \omega & -c \end{pmatrix}$; the characteristic polynomial $P(X)=(X+\delta)(X+c)(X+\mu + \beta_0 V^{is}) - \delta c (X+\mu)$ has the following properties: 
$P(-\infty)< 0$, $P(-\delta - c)=  \delta c \beta_0 V^{is} > 0$, thus it has a real root which is smaller than  $-\delta -c$. The product of all roots is $\delta c \beta_0 V^{is} >0$ and the sum of all roots is $-\delta - c - \mu - \beta_0 V^{is} < - \delta -c$, thus the other two roots have negative real part. Therefore the equilibrium is stable. 
%Any trajectory is converging to this stable equilibrium.

It is important to note that the viral load $V^{is}$ is the viral load that the infection will cause in a completely immunodeficient individual. We expect $V^{is}$  to be significantly high, see in section \ref{sec:modelvirusnoade} for details.

\subsection{Model: virus and immune response but no enhancement}
\label{sec:modelvirusnoade}
In this section we consider the model \eqref{eq:fullmodelT}-\eqref{eq:fullmodelA} with no ADE i.e.,
$\beta(A)=\beta_0$ that is $\beta_1=0$.
This model is similar to other in the literature (see for instance \cite[eq. (2.9) page 29]{wodarz_killer_2007_book} who consider also the cytotoxic effect of the immune response on the infected cells; however they do not consider virus destruction by antibodies. In particular there virus load is constant.
Another similar model is \cite[eqns. (8.1)-(8.3)]{wodarz_killer_2007_book}. 
With respect to the previous section here the immune response is present. It is triggered by a threshold set at 
{
$V^t$ (see definition in \eqref{eq:defVisVstar})}. 
It is natural to suppose that the immune response threshold 
is a very small value and in any case a value smaller than the immunosupression viral load $V^{is}$ in \eqref{eq:immunosuppression_eq}. That is we can make the hypothesis that $V^{is} > V^t$ i.e.
{
assumption 	\eqref{eq:hypothesisV} holds.
}%end modif

The Jacobian matrix is: 
\begin{equation}
	J= \begin{pmatrix} 
		-(\mu+\beta_0 V^{is}) & 0 & -\beta_0 T & 0\\
		\beta_0 V & -\delta & \beta_0 T & 0
		\\ 0 & \omega & -(c+bA) & -bV \\
		0& 0&aA& aV-\sigma
	\end{pmatrix}.
\end{equation}

With these provisions, one can find analytically the 
critical points (equilibria candidates): 
\begin{enumerate}
	\item $T=T^*=\Lambda/\mu$, $V=I=A=0$, which is the high dimensional analog of equilibrium \eqref{eq:equilibrium1d}.
	However, unlike in section \ref{sec:noinfectionmodel}, this equilibrium is not stable any more (the determinant of the Jacobian matrix is negative when hypothesis \eqref{eq:condition} is satisfied.
	\item the immunosuppression equilibrium \eqref{eq:immunosuppression_eq} with $A=0$; again this equilibrium is not stable any more because the condition \eqref{eq:hypothesisV} implies that the eigenvalue $aV^{is}- \sigma$ is positive.
	\item the only critical point left is \eqref{eq:equilibriumbeta1nul} {We prove that it is stable in the following.}	Note that the equilibrium value of the antibody level is positive due to condition \eqref{eq:hypothesisV}.
	
\end{enumerate}

{
\begin{proof}[{\bf Proof of the Proposition \ref{prop:stabilitynoade}}]
The equilibrium is stable when the real parts of the eigenvalues of the Jacobian matrix are negative. This is the same as saying that the roots of the polynomials $P_0(X)=det(X \cdot I - J)$
have negative real parts (here $I$ is the identity matrix). Such a polynomial is called stable and, if we write $P_0(X) = \gamma_4 X^4 +\gamma_3 X^3 +\gamma_2 X^2 + \gamma_1 X^1 +
\gamma_0$ then, following the Routh-Hurwitz criterion 
\cite{routh1877treatise}, 
\cite[p. 1076]{routh_hurwitz}, the stability holds true when
\begin{equation}
\gamma_k >0, k=0,1,2,3,4
\label{eq:routh4_condition_positive_coeffs}
\end{equation}
and
\begin{equation}
\gamma_1 \gamma_2 \gamma_3 >  \gamma_4 \gamma_1^2 + \gamma_3^2 \gamma_0.
\label{eq:main_routh4_condition}
\end{equation}
Unfortunately, checking in general these conditions is very difficult because the expressions involved are highly non-linear in the original parameters of the model ($a$,$b$, $c$, $\sigma$, etc.). We therefore need to exploit to the full extent the specific setting of our model. To this end we will make the following change of variables:
\begin{equation}
\zeta= \delta - \mu > 0, \ \ 
%\xi = R_0 - 1 > 0 \ \ 
%w = \xi-\frac{\sigma \beta_0}{\mu a}
w = R_0 - 1-\frac{\sigma \beta_0}{\mu a}>0.
%w = \left( R_0 - 1 \right) - \frac{\beta_0 \sigma}{\mu a}.
\end{equation}
After replacing all new variables and direct computations, we obtain:

$\gamma_0= c \delta \mu \sigma w$,
$\gamma_1= \frac{c \sigma \left(a^{2} \delta \mu w + a^{2} \mu^{2} w + a \beta_{0} \delta \mu w + a \beta_{0} \delta \mu + a \beta_{0} \mu \sigma w + \beta_{0}^{2} \delta \sigma\right)}{a \left(a \mu + \beta_{0} \sigma\right)}$,

$\gamma_2= \frac{a^{2} c \mu^{2} w + a^{2} c \mu^{2} + a^{2} c \mu \sigma w + a^{2} \delta \mu^{2} + a \beta_{0} c \mu \sigma w + 2 a \beta_{0} c \mu \sigma + 2 a \beta_{0} \delta \mu \sigma + \beta_{0}^{2} c \sigma^{2} + \beta_{0}^{2} \delta \sigma^{2}}{a \left(a \mu + \beta_{0} \sigma\right)}$, 

$\gamma_3=  \frac{a^{2} \delta \mu + a^{2} \mu^{2} + a \beta_{0} \delta \sigma + 2 a \beta_{0} \mu \sigma + a c \left(a \mu \left(w + 1\right) + \beta_{0} \sigma\right) + \beta_{0}^{2} \sigma^{2}}{a \left(a \mu + \beta_{0} \sigma\right)}$, $\gamma_4= 1$.

\noindent
Since all parameters involved are positive we obtain that the condition \eqref{eq:routh4_condition_positive_coeffs} is satisfied. To check the remaining condition
\eqref{eq:main_routh4_condition} we obtain
\begin{equation}
\gamma_1 \gamma_2 \gamma_3 -  \gamma_4 \gamma_1^2 + \gamma_3^2 \gamma_0 = 
\frac{a c}{a^{3} \left(a \mu + \beta_{0} \sigma\right)^{3}}
Q_0(w,a,c,\beta_{0},\mu,\sigma,\zeta),
\end{equation}
where the multi-variable polynomial $Q_0$ is seen, after long but
straightforward computations, to be equal to~:

\noindent
$Q_0(w,a,c,\beta_{0},\mu,\sigma,\zeta)=
w^{3}a^{6}c^{2}\mu^{5} + 2 w^{3}a^{6}c^{2}\mu^{4}\sigma + w^{3}a^{6}c^{2}\mu^{3}\sigma\zeta + w^{3}a^{5}c^{2}\beta_{0}\mu^{5} + 3 w^{3}a^{5}c^{2}\beta_{0}\mu^{4}\sigma + w^{3}a^{5}c^{2}\beta_{0}\mu^{4}\zeta + w^{3}a^{5}c^{2}\beta_{0}\mu^{3}\sigma^{2} + w^{3}a^{5}c^{2}\beta_{0}\mu^{3}\sigma\zeta + w^{3}a^{4}c^{2}\beta_{0}^{2}\mu^{4}\sigma + w^{3}a^{4}c^{2}\beta_{0}^{2}\mu^{3}\sigma^{2} + w^{3}a^{4}c^{2}\beta_{0}^{2}\mu^{3}\sigma\zeta + 2 w^{2}a^{6}c^{2}\mu^{5} + 2 w^{2}a^{6}c^{2}\mu^{4}\sigma + w^{2}a^{6}c^{2}\mu^{3}\sigma\zeta + 2 w^{2}a^{6}c\mu^{6} + w^{2}a^{6}c\mu^{5}\zeta + 3 w^{2}a^{5}c^{2}\beta_{0}\mu^{5} + 8 w^{2}a^{5}c^{2}\beta_{0}\mu^{4}\sigma + 3 w^{2}a^{5}c^{2}\beta_{0}\mu^{4}\zeta + 3 w^{2}a^{5}c^{2}\beta_{0}\mu^{3}\sigma^{2} + 2 w^{2}a^{5}c^{2}\beta_{0}\mu^{3}\sigma\zeta + w^{2}a^{5}c^{2}\beta_{0}\mu^{2}\sigma^{2}\zeta + 3 w^{2}a^{5}c\beta_{0}\mu^{6} + 5 w^{2}a^{5}c\beta_{0}\mu^{5}\sigma + 5 w^{2}a^{5}c\beta_{0}\mu^{5}\zeta + 2 w^{2}a^{5}c\beta_{0}\mu^{4}\zeta^{2} - w^{2}a^{5}c\beta_{0}\mu^{3}\sigma\zeta^{2} + 6 w^{2}a^{4}c^{2}\beta_{0}^{2}\mu^{4}\sigma + 8 w^{2}a^{4}c^{2}\beta_{0}^{2}\mu^{3}\sigma^{2} + 6 w^{2}a^{4}c^{2}\beta_{0}^{2}\mu^{3}\sigma\zeta + w^{2}a^{4}c^{2}\beta_{0}^{2}\mu^{2}\sigma^{3} + 2 w^{2}a^{4}c^{2}\beta_{0}^{2}\mu^{2}\sigma^{2}\zeta + 6 w^{2}a^{4}c\beta_{0}^{2}\mu^{5}\sigma + 6 w^{2}a^{4}c\beta_{0}^{2}\mu^{4}\sigma^{2} + 9 w^{2}a^{4}c\beta_{0}^{2}\mu^{4}\sigma\zeta - w^{2}a^{4}c\beta_{0}^{2}\mu^{3}\sigma^{2}\zeta + 3 w^{2}a^{4}c\beta_{0}^{2}\mu^{3}\sigma\zeta^{2} - w^{2}a^{4}c\beta_{0}^{2}\mu^{2}\sigma^{2}\zeta^{2} + 3 w^{2}a^{3}c^{2}\beta_{0}^{3}\mu^{3}\sigma^{2} + 2 w^{2}a^{3}c^{2}\beta_{0}^{3}\mu^{2}\sigma^{3} + 3 w^{2}a^{3}c^{2}\beta_{0}^{3}\mu^{2}\sigma^{2}\zeta + 4 w^{2}a^{3}c\beta_{0}^{3}\mu^{4}\sigma^{2} + 4 w^{2}a^{3}c\beta_{0}^{3}\mu^{3}\sigma^{3} + 5 w^{2}a^{3}c\beta_{0}^{3}\mu^{3}\sigma^{2}\zeta + w^{2}a^{3}c\beta_{0}^{3}\mu^{2}\sigma^{2}\zeta^{2} + w^{2}a^{2}c\beta_{0}^{4}\mu^{3}\sigma^{3} + w^{2}a^{2}c\beta_{0}^{4}\mu^{2}\sigma^{4} + w^{2}a^{2}c\beta_{0}^{4}\mu^{2}\sigma^{3}\zeta + wa^{6}c^{2}\mu^{5} + 2 wa^{6}c\mu^{6} + wa^{6}c\mu^{5}\zeta + 3 wa^{5}c^{2}\beta_{0}\mu^{5} + 5 wa^{5}c^{2}\beta_{0}\mu^{4}\sigma + 3 wa^{5}c^{2}\beta_{0}\mu^{4}\zeta + wa^{5}c^{2}\beta_{0}\mu^{3}\sigma\zeta + 6 wa^{5}c\beta_{0}\mu^{6} + 7 wa^{5}c\beta_{0}\mu^{5}\sigma + 10 wa^{5}c\beta_{0}\mu^{5}\zeta + wa^{5}c\beta_{0}\mu^{4}\sigma\zeta + 4 wa^{5}c\beta_{0}\mu^{4}\zeta^{2} - wa^{5}c\beta_{0}\mu^{3}\sigma\zeta^{2} + 2 wa^{5}\beta_{0}\mu^{7} + 5 wa^{5}\beta_{0}\mu^{6}\zeta + 4 wa^{5}\beta_{0}\mu^{5}\zeta^{2} + wa^{5}\beta_{0}\mu^{4}\zeta^{3} + 9 wa^{4}c^{2}\beta_{0}^{2}\mu^{4}\sigma + 8 wa^{4}c^{2}\beta_{0}^{2}\mu^{3}\sigma^{2} + 9 wa^{4}c^{2}\beta_{0}^{2}\mu^{3}\sigma\zeta + 2 wa^{4}c^{2}\beta_{0}^{2}\mu^{2}\sigma^{2}\zeta + 18 wa^{4}c\beta_{0}^{2}\mu^{5}\sigma + 11 wa^{4}c\beta_{0}^{2}\mu^{4}\sigma^{2} + 28 wa^{4}c\beta_{0}^{2}\mu^{4}\sigma\zeta - wa^{4}c\beta_{0}^{2}\mu^{3}\sigma^{2}\zeta + 10 wa^{4}c\beta_{0}^{2}\mu^{3}\sigma\zeta^{2} - 2 wa^{4}c\beta_{0}^{2}\mu^{2}\sigma^{2}\zeta^{2} + 7 wa^{4}\beta_{0}^{2}\mu^{6}\sigma + 17 wa^{4}\beta_{0}^{2}\mu^{5}\sigma\zeta + 13 wa^{4}\beta_{0}^{2}\mu^{4}\sigma\zeta^{2} + 3 wa^{4}\beta_{0}^{2}\mu^{3}\sigma\zeta^{3} + 9 wa^{3}c^{2}\beta_{0}^{3}\mu^{3}\sigma^{2} + 5 wa^{3}c^{2}\beta_{0}^{3}\mu^{2}\sigma^{3} + 9 wa^{3}c^{2}\beta_{0}^{3}\mu^{2}\sigma^{2}\zeta + wa^{3}c^{2}\beta_{0}^{3}\mu\sigma^{3}\zeta + 20 wa^{3}c\beta_{0}^{3}\mu^{4}\sigma^{2} + 10 wa^{3}c\beta_{0}^{3}\mu^{3}\sigma^{3} + 28 wa^{3}c\beta_{0}^{3}\mu^{3}\sigma^{2}\zeta - wa^{3}c\beta_{0}^{3}\mu^{2}\sigma^{3}\zeta + 8 wa^{3}c\beta_{0}^{3}\mu^{2}\sigma^{2}\zeta^{2} - wa^{3}c\beta_{0}^{3}\mu\sigma^{3}\zeta^{2} + 9 wa^{3}\beta_{0}^{3}\mu^{5}\sigma^{2} + 21 wa^{3}\beta_{0}^{3}\mu^{4}\sigma^{2}\zeta + 15 wa^{3}\beta_{0}^{3}\mu^{3}\sigma^{2}\zeta^{2} + 3 wa^{3}\beta_{0}^{3}\mu^{2}\sigma^{2}\zeta^{3} + 3 wa^{2}c^{2}\beta_{0}^{4}\mu^{2}\sigma^{3} + wa^{2}c^{2}\beta_{0}^{4}\mu\sigma^{4} + 3 wa^{2}c^{2}\beta_{0}^{4}\mu\sigma^{3}\zeta + 10 wa^{2}c\beta_{0}^{4}\mu^{3}\sigma^{3} + 5 wa^{2}c\beta_{0}^{4}\mu^{2}\sigma^{4} + 12 wa^{2}c\beta_{0}^{4}\mu^{2}\sigma^{3}\zeta + 2 wa^{2}c\beta_{0}^{4}\mu\sigma^{3}\zeta^{2} + 5 wa^{2}\beta_{0}^{4}\mu^{4}\sigma^{3} + 11 wa^{2}\beta_{0}^{4}\mu^{3}\sigma^{3}\zeta + 7 wa^{2}\beta_{0}^{4}\mu^{2}\sigma^{3}\zeta^{2} + wa^{2}\beta_{0}^{4}\mu\sigma^{3}\zeta^{3} + 2 wac\beta_{0}^{5}\mu^{2}\sigma^{4} + wac\beta_{0}^{5}\mu\sigma^{5} + 2 wac\beta_{0}^{5}\mu\sigma^{4}\zeta + wa\beta_{0}^{5}\mu^{3}\sigma^{4} + 2 wa\beta_{0}^{5}\mu^{2}\sigma^{4}\zeta + wa\beta_{0}^{5}\mu\sigma^{4}\zeta^{2} + a^{5}c^{2}\beta_{0}\mu^{5} + a^{5}c^{2}\beta_{0}\mu^{4}\zeta + 3 a^{5}c\beta_{0}\mu^{6} + 5 a^{5}c\beta_{0}\mu^{5}\zeta + 2 a^{5}c\beta_{0}\mu^{4}\zeta^{2} + 2 a^{5}\beta_{0}\mu^{7} + 5 a^{5}\beta_{0}\mu^{6}\zeta + 4 a^{5}\beta_{0}\mu^{5}\zeta^{2} + a^{5}\beta_{0}\mu^{4}\zeta^{3} + 4 a^{4}c^{2}\beta_{0}^{2}\mu^{4}\sigma + 4 a^{4}c^{2}\beta_{0}^{2}\mu^{3}\sigma\zeta + 12 a^{4}c\beta_{0}^{2}\mu^{5}\sigma + 19 a^{4}c\beta_{0}^{2}\mu^{4}\sigma\zeta + 7 a^{4}c\beta_{0}^{2}\mu^{3}\sigma\zeta^{2} + 9 a^{4}\beta_{0}^{2}\mu^{6}\sigma + 22 a^{4}\beta_{0}^{2}\mu^{5}\sigma\zeta + 17 a^{4}\beta_{0}^{2}\mu^{4}\sigma\zeta^{2} + 4 a^{4}\beta_{0}^{2}\mu^{3}\sigma\zeta^{3} + 6 a^{3}c^{2}\beta_{0}^{3}\mu^{3}\sigma^{2} + 6 a^{3}c^{2}\beta_{0}^{3}\mu^{2}\sigma^{2}\zeta + 19 a^{3}c\beta_{0}^{3}\mu^{4}\sigma^{2} + 28 a^{3}c\beta_{0}^{3}\mu^{3}\sigma^{2}\zeta + 9 a^{3}c\beta_{0}^{3}\mu^{2}\sigma^{2}\zeta^{2} + 16 a^{3}\beta_{0}^{3}\mu^{5}\sigma^{2} + 38 a^{3}\beta_{0}^{3}\mu^{4}\sigma^{2}\zeta + 28 a^{3}\beta_{0}^{3}\mu^{3}\sigma^{2}\zeta^{2} + 6 a^{3}\beta_{0}^{3}\mu^{2}\sigma^{2}\zeta^{3} + 4 a^{2}c^{2}\beta_{0}^{4}\mu^{2}\sigma^{3} + 4 a^{2}c^{2}\beta_{0}^{4}\mu\sigma^{3}\zeta + 15 a^{2}c\beta_{0}^{4}\mu^{3}\sigma^{3} + 20 a^{2}c\beta_{0}^{4}\mu^{2}\sigma^{3}\zeta + 5 a^{2}c\beta_{0}^{4}\mu\sigma^{3}\zeta^{2} + 14 a^{2}\beta_{0}^{4}\mu^{4}\sigma^{3} + 32 a^{2}\beta_{0}^{4}\mu^{3}\sigma^{3}\zeta + 22 a^{2}\beta_{0}^{4}\mu^{2}\sigma^{3}\zeta^{2} + 4 a^{2}\beta_{0}^{4}\mu\sigma^{3}\zeta^{3} + ac^{2}\beta_{0}^{5}\mu\sigma^{4} + ac^{2}\beta_{0}^{5}\sigma^{4}\zeta + 6 ac\beta_{0}^{5}\mu^{2}\sigma^{4} + 7 ac\beta_{0}^{5}\mu\sigma^{4}\zeta + ac\beta_{0}^{5}\sigma^{4}\zeta^{2} + 6 a\beta_{0}^{5}\mu^{3}\sigma^{4} + 13 a\beta_{0}^{5}\mu^{2}\sigma^{4}\zeta + 8 a\beta_{0}^{5}\mu\sigma^{4}\zeta^{2} + a\beta_{0}^{5}\sigma^{4}\zeta^{3} + c\beta_{0}^{6}\mu\sigma^{5} + c\beta_{0}^{6}\sigma^{5}\zeta + \beta_{0}^{6}\mu^{2}\sigma^{5} + 2 \beta_{0}^{6}\mu\sigma^{5}\zeta + \beta_{0}^{6}\sigma^{5}\zeta^{2}$.

Most of the monomials in $Q_0$ have positive coefficients, except 
the following ones:
$-w^{2}a^{5}c\beta_{0}\mu^{3}\sigma\zeta^{2}$, $- w^{2}a^{4}c\beta_{0}^{2}\mu^{3}\sigma^{2}\zeta$,
$- w^{2}a^{4}c\beta_{0}^{2}\mu^{2}\sigma^{2}\zeta^{2}$, 
$- wa^{5}c\beta_{0}\mu^{3}\sigma\zeta^{2}$, 
$- wa^{4}c\beta_{0}^{2}\mu^{3}\sigma^{2}\zeta$, 
$- 2 wa^{4}c\beta_{0}^{2}\mu^{2}\sigma^{2}\zeta^{2}$,
$- wa^{3}c\beta_{0}^{3}\mu^{2}\sigma^{3}\zeta$, 
$- wa^{3}c\beta_{0}^{3}\mu\sigma^{3}\zeta^{2}$.
However, in all cases we can come up with two terms that render the total sum positive. 
For instance
the term $-w^{2}a^{5}c\beta_{0}\mu^{3}\sigma\zeta^{2}$
(term $27$ of the polynomial)
 is negative but, when
we combine it with the terms 
$w^{3}a^{6}c^2\mu^{3}\sigma\zeta /2 $ (half of the third term) 
%2: Out[94]: array([3, 6, 2, 0, 3, 1, 1])
and 
% 78: Out[95]: array([1, 4, 0, 2, 3, 1, 3])
$w a^{4}\beta_{0}^2\mu^{3}\sigma\zeta^{3}/2$ (half of the term $79$),
both appearing with positive coefficients, we obtain a positive number
$w^{3}a^{6}c^2\mu^{3}\sigma\zeta/2 -w^{2}a^{5}c\beta_{0}\mu^{3}\sigma\zeta^{2}+
w a^{4}\beta_{0}^2\mu^{3}\sigma\zeta^{3}/2 
= \frac{wa\mu^{3}\sigma\zeta}{2}( w ac -\beta_0 \zeta )^2
\ge 0$.

The interested reader can check that in the same way that:

- the  $36$-th monomial compensate with monomials $7$ and $92$;

- the  $38$-th monomial compensate with monomials $3$ and $104$;

- the $61$-th monomial compensate with monomials $14$ and $128$;

- the $73$-th monomial compensate with monomials $13$ and $145$;

- the $75$-th monomial compensate with monomials $14$ and $146$;

- the $87$-th monomial compensate with monomials $20$ and $154$;

- the $89$-th monomial compensate with monomials $22$ and $155$.

This allows to state that $Q_0>0$ which concludes the proof. 
%A symbolic computation program
%that allows to check this proof is provided (using the "Sympy" Python library \cite{sympy}).
\end{proof}

\subsection{Full model: virus, immune system and ADE}
\label{sec:modelvirus_full}

\begin{proof}[Proof of the Proposition \ref{prop:stability_full}]
We consider the model \eqref{eq:fullmodelT}-\eqref{eq:fullmodelA} with $\beta(A)=\beta_0+ \beta_1 A $ ($\beta_1>0$).
The analysis of this dynamics is more involved. The first two equilibria, having $A=0$ are the complete analogues of the 
equilibria seen in previous sections and have no dynamical interest. Since $A=0$ the parameter $\beta_1$ that multiplies $A$ has no impact and the proof of the instability of the trivial equilibrium and immunosupression equilibrium follow exactly the same arguments as before.

To find the third equilibrium, note that after immediate computations we find that  the antibody level is solution of the second order equation \eqref{eq:equationAordre2}. Such an equation has two solutions but exactly one is positive because the product of roots is negative; thus only a single point is an admissible  equilibrium, namely the positive solution of \eqref{eq:equationAordre2} (with respect to the unknown $A$); setting to zero all
derivatives we obtain the other values as in 
\eqref{eq:equilibriumbeta1positive}.

To prove the properties of this equilibrium we start with the point 
\ref{item:prop_counter_example} of the proposition; consider the values 
$a=\sigma=c=b=\omega=1$, $\mu=1.e-3$, $\delta=2$, $\Lambda=4$,
$\beta_0=0.0011$ and $\beta_1=0.01188$; all hypotheses are satisfied and the numerical 
values of the equilibrium are $T=333.33, I=1.83, V=1, A=0.83$ while the eigenvalues are $-3.45$, $0.50$, $0.01$ and $-0.90$. Since some eigenvalues are real and positive the equilibrium is not stable for this set of parameters. This completes the proof for this point.
In practice the evolution oscillates indefinitely between a state with high $T$ value and one with very low $T$ value.

Note that the point \ref{item:prop_beta1_small} of the conclusion is just a consequence of the continuity and the proposition \ref{prop:stabilitynoade}, because both the equilibrium and the coefficients of the polynomial $P(X) = \det(X \times Id - J)$ evaluated at the equilibrium depend smoothly on $\beta_1$. 
Since we proved  that 
\eqref{eq:routh4_condition_positive_coeffs}
and \eqref{eq:main_routh4_condition}
are true for $\beta_1=0$ by continuity the terms in the two conditions will remain strictly positive for $\beta_1$ small enough and by the Routh-Hurwitz criterion the equilibrium will be stable.

The only point remaining to be proved is \ref{item:prop_beta1_large}.
Note that when $\beta_1 \to \infty$ the positive root $A^f$ of the equation \eqref{eq:equationAordre2}  
converges to some quantity $A^\infty$, and $\beta(A^f) \to \infty$; moreover, 
we obtain from the definition of $T^f$ that 
$\lim_{\beta_1 \to \infty} T^f = 0$ and 
$\lim_{\beta_1 \to \infty} (\beta(A^f) T^f)=\lim_{\beta_1 \to \infty} (\beta_1 T^f)=\frac{\delta (c+bA^\infty)}{\omega A^\infty}$.
% will converge to $\beta^\infty = \frac{\delta c \mu}{\omega \Lambda - \delta c V^t}$.
Consider the Jacobian matrix: 
{\scriptsize
\begin{equation}
	J= \begin{pmatrix} 
		-\beta(A)V-\mu & 0 & -\beta(A) T & -\beta_1 T V\\
		\beta(A) V & -\delta & \beta(A) T & \beta_1 T V
		\\ 0 & \omega & -(c+bA) & -bV \\
		0& 0&aA& aV-\sigma
	\end{pmatrix}.
\end{equation}
}
Let us compute $P(X)=\det(X Id - J)=\det(J-X Id)$:
{\scriptsize
\begin{eqnarray}
	& \ & 
	P(X)= \left|\begin{matrix} 
		-\beta(A)V-\mu-X & 0 & -\beta(A) T & -\beta_1 T V\\
		\beta(A) V & -\delta-X & \beta(A) T & \beta_1 T V
		\\ 0 & \omega & -X-(c+bA) & -bV \\
		0& 0&aA& -X
	\end{matrix}\right|
\nonumber \\ & \ & = 
\left|\begin{matrix} 
	-\mu-X & -\delta-X & 0 & 0\\
	\beta(A) V & -\delta-X & \beta(A) T & \beta_1 T V
	\\ 0 & \omega & -X-(c+bA) & -bV \\
	0& 0&aA& -X
\end{matrix}\right|
\nonumber \\ & \ & =
-(\mu+X) \left|\begin{matrix} 
 -\delta-X & \beta(A) T & \beta_1 T V
	\\  \omega & -X-(c+bA) & -bV \\
 0&aA& -X
\end{matrix}\right|
- \beta(A) V 
\left|\begin{matrix} 
 -\delta-X & 0 & 0
	\\  \omega & -X-(c+bA) & -bV \\
	 0&aA& -X
\end{matrix}\right|.\nonumber
\end{eqnarray}
}

Thus the polynomial $P(X)=\det(X Id - J)$ can be written, to first order in $\beta_1$, as 
$P(X) = R(X) + \beta(A) V(X + \delta)(X^2 + X(c+bA)+ aAbV)$, 
where $R(X)$ is a fourth order polynomial with leading term $X^4$ and coefficients independent
of $\beta_1$. 
Note that $(X + \delta)(X^2 + X(c+bA)+ aAbV)$ is a stable polynomial. 
To finish the proof we invoke Lemma~\ref{lemma:stability3} below for 
$\psi=\beta(A)V$.

\end{proof}
\begin{lemma}
	Let $Z_3 = \phi_3 X^3 + \phi_2 X^2 + \phi_1 X + \phi_0$ be a stable polynomial of order $3$ with $\phi_3>0$ and 
	$Z_4 = \varphi_4 X^4+ \varphi_3 X^3 + \varphi_2 X^2 + \varphi_1 X + \varphi_0$ a polynomial of order four with $\varphi_4>0$. Then, for $\psi$ large enough
	the polynomial $Z_4(X) + \psi Z_3(X)$ is stable.
	\label{lemma:stability3}
\end{lemma}
\begin{proof}
	Since $\phi_3>0$
	using the reciprocal of the Routh-Hurwitz criterion 
	all coefficients $\phi_k$ are strictly positive and 
	$\phi_1 \phi_2 > \phi_0 \phi_3$. For $\psi$ large enough this allows to check the Routh-Hurwitz
	criterion for the fourth order polynomial 
	$Z_4(X) + \psi Z_3(X)$: the coefficients will be positive and the last remaining condition 
	is 
	$(\varphi_1 + \psi \phi_1) (\varphi_2 + \psi \phi_2)
	(\varphi_3 + \psi \phi_3) > (\varphi_0 + \psi \phi_0) (\varphi_3 + \psi \phi_3)^2 
	+(\varphi_1 + \psi \phi_1) \varphi_4 ^2$, which is satisfied for $\psi$ large enough (leading term 
	$(\phi_1 \phi_2- \phi_0 \phi_3)\phi_3$
	is positive).
\end{proof}

\section{Extended model including a latent phase} \label{sec:extended_latent_model}

We present here a version of the main model \eqref{eq:fullmodelT}-\eqref{eq:fullmodelbeta}
extended to take into account a latent phase of the cells. The interest of such a model is to
 give a finer description of all states of the attacked cells; this comes however at the price of requiring several more parameters (including the transition rate $\eta>0$ from the latent to infected, virus-producing, cells). In practice the choice of the model depends on the outcomes of interest and available data to fit. In our case the data to fit was relatively scarce thus we kept the restricted model \eqref{eq:fullmodelT}-\eqref{eq:fullmodelbeta} for the numerical simulations. Denoting $L$ the 
 number of latent infected cells (i.e., cells already infected but not yet producing viruses) we can write this model as:

\begin{eqnarray}
	& \ & dT/dt = \Lambda - \mu T -\beta(A) V T	\label{eq:latentmodelT}		\\ 
	& \ & dL/dt = \beta(A)V T - \eta L- \mu L	\label{eq:latentmodelL}		\\ 
	& \ & dI/dt =  \eta L - \delta I  \label{eq:latentmodelI}\\
	 & \ & 
	dV / dt = \omega I - c V - b A V  \label{eq:latentmodelV}
	\\ & \ & 
	dA/dt = a V A - \sigma A \label{eq:latentmodelA}
	\\ & \ & 
	\beta(A)=\beta_0+\beta_1 A. \label{eq:latentmodelbeta}
\end{eqnarray}
The model is illustrated in figure 	\ref{fig:tliva_latent_model}.

\begin{figure}[htbp!]
	\centering
	\includegraphics[width=0.9\textwidth]{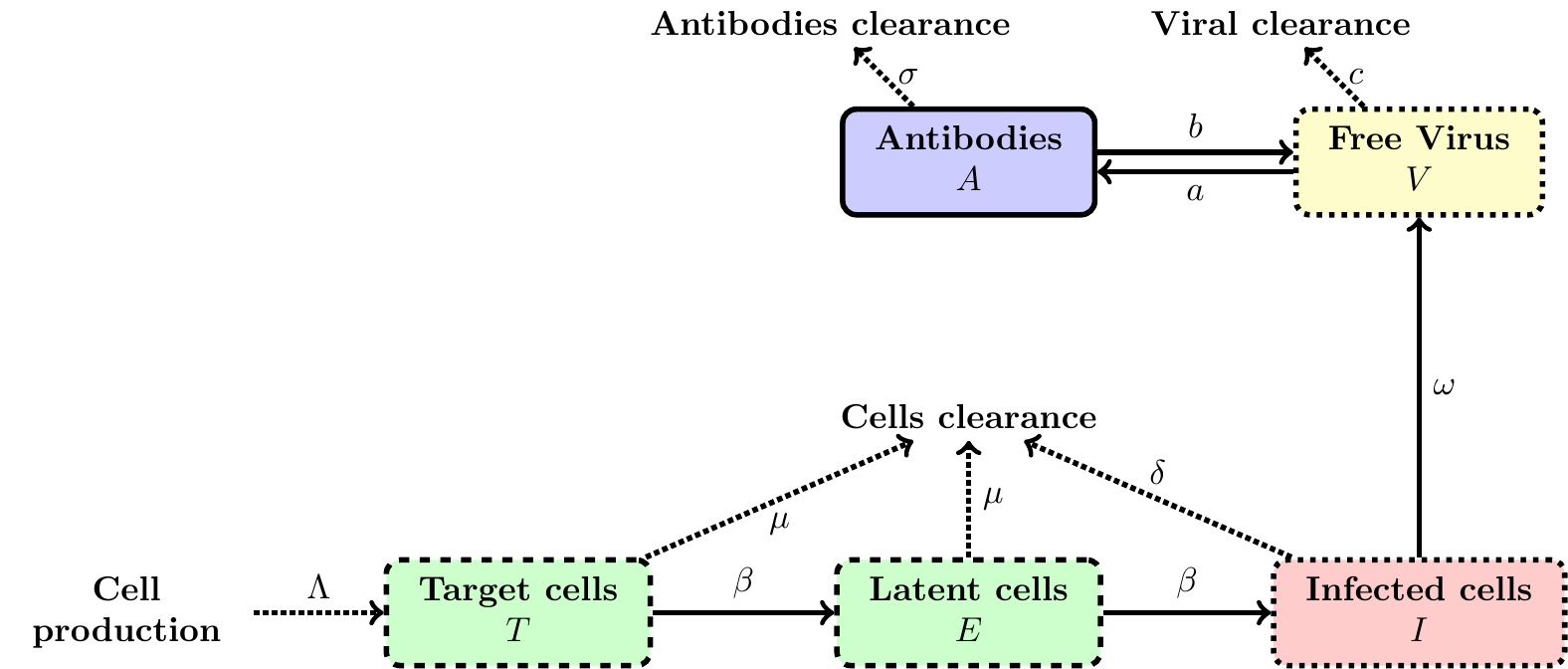}
	\caption{Graphical illustration of the flow in the model \eqref{eq:latentmodelT}-\eqref{eq:latentmodelbeta}.}
	\label{fig:tliva_latent_model}
\end{figure}
}
\end{document}